\documentclass[journal,twoside]{IEEEtran}

\usepackage{graphicx}
\usepackage{epstopdf}

\def\e{\mathop{\rm e}\nolimits}
\def\sgn{\mathop{\rm sgn}\nolimits}
\def\Re{\mathop{\rm Re}\nolimits}
\def\Im{\mathop{\rm Im}\nolimits}

\def\erfi{\mathop{\rm erfi}\nolimits}

\def\Xint#1{\mathchoice
   {\XXint\displaystyle\textstyle{#1}}%
   {\XXint\textstyle\scriptstyle{#1}}%
   {\XXint\scriptstyle\scriptscriptstyle{#1}}%
   {\XXint\scriptscriptstyle\scriptscriptstyle{#1}}%
   \!\int}
\def\XXint#1#2#3{{\setbox0=\hbox{$#1{#2#3}{\int}$}
     \vcenter{\hbox{$#2#3$}}\kern-.5\wd0}}

\def\dashint{\Xint-}

\font\bb=msbm10 scaled \magstep1 
\def\R{\hbox{\bb R}} 
 
\def\Z{\hbox{\bb Z}} 

\newtheorem{theorem}{Theorem}{}

\begin{document}
%
\title{Time Delay Extraction from Frequency Domain Data Using Causal Fourier Continuations 
for High-Speed Interconnects
}



\markboth
{Lyudmyla~L.~Barannyk \MakeLowercase{\textit{et al.}}: Time Delay Extraction from Frequency Domain Data Using Causal Fourier Continuations}{Electronics}%


\author{Lyudmyla~L.~Barannyk${}^{1}$,
Hung~H.~Tran${}^2$, 
        Aicha~Elshabini${}^2$,
        and~Fred~D.~Barlow${}^3$
\thanks{Lyudmyla~L.~Barannyk is 
with the Department of Mathematics, University of Idaho, Moscow, Idaho 83844, USA (e-mail: barannyk@uidaho.edu).}
\thanks{Hung~H.~Tran and Aicha~Elshabini are with the Department of Electrical \& Computer Engineering, University of Idaho, Moscow, Idaho 83844, USA (e-mails: tran4105@vandals.uidaho.edu; elshabini@uidaho.edu).}
\thanks{Fred Barlow is with the College of Engineering,  University of Alaska Anchorage, Anchorage, Alaska 99508, USA (e-mail: fdbarlow@uaa.alaska.edu).}
}


%


\maketitle

\begin{abstract}

We present a new method for time delay estimation using band limited frequency domain data representing the port responses of interconnect structures. The approach is based on the recently developed by the authors spectrally accurate method for causality characterization that employs SVD-based causal Fourier continuations. The time delay extraction  is constructed by incorporating a linearly varying phase factor to the system of equations that determines Fourier coefficients. 
The method is capable of determining time delay using data affected by noise or approximation errors that come from measurements or numerical simulations. It can also be employed when only a limited number of frequency responses is available. The technique can be extended to multi-port and mixed mode networks. Several analytical and simulated examples are used to demonstrate the accuracy and strength of the proposed technique.

%

\end{abstract}



\begin{IEEEkeywords}
time delay, delay estimation, causality, dispersion relations, singular value decomposition, SVD-based causal Fourier continuation, high speed interconnects.
\end{IEEEkeywords}

%
\IEEEpeerreviewmaketitle

\section{Introduction}

%

Identification and extraction of time delay is an important research problem in signal processing and has applications in many fields including radar \cite{Li_2004}, sonar \cite{Quazi_1981, Chen_Benesty_Huang_2006}, ultrasonics \cite{De_Marchi_Marzani_Speciale_2009}, microwave imaging \cite{Lim_Nhung_Li_Thang_2008}, geophysics \cite{Kepko_Kivelson_1999}, seismology \cite{Wilcock_2004, Mercerat_Nolet_2013}, wireless communications \cite{Vanderveen_Vanderveen_Paulraj_1998} as well as  modeling of passive structures in electronic systems, in particular, transmission line modeling \cite{Gustavsen_2004, De_Tommasi_Gustavsen_2006},
transient simulation of interconnects \cite{Lalgudi_Engin_Casinovi_Swaminathan_2008} and co-simulation of passive structures with active devices in a time domain using SPICE. 
%
%
%
%
%
Passive structures in electronic systems have been traditionally analyzed in the frequency domain, while transient simulations are performed in the time domain using suitable  models that accurately capture  the relevant electromagnetic phenomena. The models are obtained from either direct measurements or electromagnetic simulations. Interconnect models are typically approximated by rational transfer functions using the vector fitting algorithm in various implementations \cite{Gustavsen_Semlyen_1999, Grivet_Talocia_2003, Gustavsen_Semlyen_2004, Deschrijver_Dhaene_2005, Grivet_Talocia_Bandinu_2006, Deschrijver_Haegeman_Dhaene_2007, Chinea_Triverio_Grivet_Talocia_2010}, which is the standard macromodeling approach.
As clock frequencies increase, the size of passive structures becomes of the same order as the signal wavelength at the operating frequency, which causes the distributed effects such as time delay to play a significant role in the time domain simulations. For this reason, time delay has to be included in  macromodeling, in particular, when causality is analyzed. The connection between causality and time delay is in the fact that time delays can pull a non-causal signal  into the causal region or vice versa pull a causal signal into the non-causal region, while causality, in turn, can be expressed in terms of the Hilbert transform \cite{Papoulis_1977, Weideman_1995, Nussenzveig_1972}. Several approaches can be used to extract delays in the frequency domain, for example, using the Hilbert transform \cite{Grennberg_Sandell_1994, Tsuchiya_Miki_2004, Knockaert_Dhaene_2008},
 the minimum phase all-pass decomposition \cite{Mandrekar_Swaminathan_2005, Mandrekar_Srinivasan_Engin_Swminathan_2007, Lalgudi_Engin_Casinovi_Swaminathan_2008}, 
 incorporating an optimal time delay into the vector fitting algorithm  \cite{Gustavsen_2004, De_Tommasi_Gustavsen_2006}, 
employing a modified Lie approximation to develop a passive and compact macromodel \cite{Nakhla_Dounavis_Achar_Nakhla_2005},
using a Gabor transform to develop delayed rational function macromodels for long interconnects \cite{Grivet_Talocia_2006, Chinea_Triverio_Grivet_Talocia_2010} or conducting a probabilistic analysis of the cepstrum in the presence of noise \cite{Hassab_Boucher_1976}. In the time domain, delayed rational functions \cite{Charest_Saraswat_Nakhla_Achar_Soveiko_2007, Charest_Nakhla_Achar_Saraswat_Soveiko_Erdin_2010} can be employed to extract delays. In this paper, a novel approach is proposed in which time delay is determined in the frequency domain using a causality argument. Causality is verified using the SVD-based causal Fourier continuation method developed by the authors \cite{Barannyk_Aboutaleb_Elshabini_Barlow_IEEE, Barannyk_Aboutaleb_Elshabini_Barlow_IMAPS2014}, while the time delay presence is incorporated by a linearly varying phase factor to the system of equations that determines Fourier coefficients.  Preliminary results are reported in \cite{Barannyk_Tran_Nguyen_Elshabini_Barlow_EPEPS2015}.

The rest of the paper is organized as follows. Section \ref{causality} provides a  background on causality for linear time-translation invariant systems and dispersion relations. In Section \ref{Fourier_continuation}, we show main steps in the derivation of causal  Fourier continuations using truncated singular value decomposition (SVD) method that was developed to access causality.  We also provide error estimates that take into account a possible presence of noise in data. Section \ref{delay_estimation} extends the causality characterization method to develop a technique for time delay extraction. The proposed method is tested in Section \ref{numerical_examples} using several analytic and simulated examples. We also analyze the performance of the algorithm when only a limited number of frequency  responses is available and when noise/approximation errors are present in data.
In Section \ref{conclusions} we present our conclusions. The Appendix section is devoted to formulation of error bounds for the causality characterization method based on causal Fourier continuations.

%

\section{Causality of Linear Time-Invariant Systems} \label{causality}

Consider a linear and time-invariant physical system with the impulse response  ${ h}(t)$ 
subject to a time-dependent input ${ f(t)}$, to which it responds by an output ${ x(t)}$. 
%
%
%
Denote by 
\begin{equation} \label{1_3_3}
{ H}(w)=\int_{-\infty}^\infty { h}(\tau)\e^{-i w\tau}d\tau
\end{equation}
the Fourier transform of ${ h(t)}$, which is also called the transfer function.  

The system is causal if  the output cannot precede the input, i.e. if ${ f
}(t)=0$  for $t<T$, the same must be true for ${ x(t)}$. This primitive causality condition in the time domain  implies ${ h}(t)=0$,  $t<0$. Hence, domain of integration in (\ref{1_3_3}) can be reduced to $[0,\infty)$.

Assume   $H(w)\in L_2(\R)$. Then starting from Cauchy's theorem and using contour integration, one can show \cite{Nussenzveig_1972} that for any point $w$ on the real axis, $H(w)$ can be written\footnote{Please note that we use an opposite sign of the exponent in the definition of the Fourier transform than in \cite{Nussenzveig_1972}.} as
\begin{equation} \label{1_6_7}
H(w)=\frac{1}{\pi i}\dashint_{-\infty}^\infty \frac{H(w')}{w-w'}dw', \quad \mbox{real} \ w,
\end{equation}
where 
\[
\dashint_{-\infty}^\infty=P\int_{-\infty}^\infty = \lim_{\epsilon\to 0}\left(\int_{-\infty}^{w-\epsilon}+\int_{w+\epsilon}^{\infty}\right)
\]
denotes Cauchy's principal value. Separating the real and imaginary parts of (\ref{1_6_7}), we get
\begin{equation} \label{1_6_10}
\Re H(w)=\frac{1}{\pi}\dashint_{-\infty}^\infty \frac{\Im H(w')}{w-w'}dw',
\end{equation}
\begin{equation} \label{1_6_11}
\Im H(w)=-\frac{1}{\pi}\dashint_{-\infty}^\infty \frac{\Re H(w')}{w-w'}dw'.
\end{equation}
Expressions (\ref{1_6_10}) and (\ref{1_6_11}) are called the dispersion relations or Kramers-Kr\"onig  relations. They show that $\Re H$ and $\Im H$ are not independent functions, but instead they are related to each other: $\Re H$ at one frequency depends on $\Im H$ at all frequencies, and vice versa. This implies that if one of the functions $\Re H$ or $\Im H$ is square integrable and known, then the other one can be   completely determined by causality. 
Recalling the definition of the Hilbert transform,
\[
{\mathcal H}[u(w)]=\frac{1}{\pi}\dashint_{-\infty}^\infty \frac{u(w')}{w-w'} dw',
\]
we see that $\Re H$ and $\Im H$ are Hilbert transforms of each other, i.e.
\begin{equation} \label{Hilbert_pair}
\Re H(w)={\mathcal H}[\Im H(w)], \quad \Im H(w)=-{\mathcal H}[\Re H(w)].
\end{equation}
In other words, $\Re H$ or $\Im H$  form a Hilbert transform pair. Dispersion relations provide the causality condition in the frequency domain.


Evaluation of the Hilbert transform requires integration on $(-\infty,\infty)$, which can be reduced to $[0,\infty)$ by spectrum symmetry of $H(w)$ if $h(t)$ is real valued.
 In practice, only a limited number of discrete values of  $H(w)$ is   available  on $[w_{min},w_{max}]$. Thus, the domain of integration has to be truncated. This usually causes serious boundary artifacts due to the lack of out-of-band frequency responses. To reduce or even completely remove boundary artifacts,  the authors recently developed periodic polynomial \cite{Aboutaleb_Barannyk_Elshabini_Barlow_WMED13, Barannyk_Aboutaleb_Elshabini_Barlow_IMAPS} and causal Fourier continuation \cite{Barannyk_Aboutaleb_Elshabini_Barlow_IMAPS2014, Barannyk_Aboutaleb_Elshabini_Barlow_IEEE}, respectively, based methods for causality characterization. The approach was motivated by the  example $H(w)=\e^{-iaw}$, $a>0$, that is not square integrable but still satisfies the dispersion relations. The causality characterization method based on causal Fourier continuations allows one to construct highly  accurate approximations of a given transfer function on the original frequency interval $[w_{min},w_{max}]$ with the uniform error that decreases as the number of Fourier coefficients increases. The technique is applicable to both  baseband and  bandpass cases and capable of detecting very small localized causality violations. The method can also be extended to  multidimensional cases. 
 
In the next section, for completeness of presentation, we show main steps in the derivation of the causal Fourier continuation method that can be used to access causality of a given transfer function whose values are available at a discrete set of frequencies. We also provide upper bounds of reconstruction error between the given function and its causal Fourier continuation. We use these error estimates to understand  how to extract time delay when data with different resolutions are available and when data are affected by noise or other approximation errors.
 


%

\section{Causal Fourier Continuations} \label{Fourier_continuation}

Consider a transfer function $H(w)=\Re H+i\Im H$, whose $N$ discrete values are available on $[w_{min}, w_{max}]$, $w_{min}\geq 0$. For real-valued impulse response functions $h(t)$, $\Re H$ and $\Im H$ are even and odd functions, respectively. This implies that $H(w)$ has values on $[-w_{max},-w_{min}]$ by spectrum symmetry. For convenience, we rescale the frequency interval $[-w_{max},w_{max}]$ to $[-0.5,0.5]$ by the substitution $x=\frac{0.5}{w_{max}}w$, so the rescaled transfer function $H(x)$ is defined on the unit length interval with $\tilde N$ values where  $\tilde N=2N-1$ or $\tilde N=2N$ depending if $H(x)$ is available at $x=0$ or not. Both baseband and bandpass cases can be considered.

The idea of a causal Fourier continuation is to construct an accurate Fourier series approximation of $H(x)$ by allowing the Fourier series to be periodic and causal in an extended domain. The result is the Fourier continuation of $H$ that we denote by ${\mathcal C}(H)$, and it is defined by
\begin{equation}\label{E1}
{\mathcal C}(H)(x)=\sum_{k=-M+1}^{M} \alpha_k \e^{-\frac{2\pi i}{b} k x},
\end{equation}
for even number $2M$ of terms, whereas for odd number $2M+1$ of terms, the index $k$ varies from $-M$ to $M$. 
Throughout this paper, we assume that the number $M$ of Fourier coefficients is even, for simplicity. When $M$ is odd, analogous results can be formulated.
%
Here $b$ is the period of approximation. For SVD-based periodic continuations $b$ is normally chosen as twice the length of the domain on which function $H$ is given though the value $b=2$ is not necessarily optimal. The optimal value $b$ depends on a function being approximated. In practice, several values $b\in(1,4)$ may be tried to get a better reconstruction of $H(x)$ with a Fourier series.

Functions $\phi_{k}(x)=\e^{-\frac{2\pi i }{b}k x}$, $k\in\Z$, form a complete orthogonal basis in $L_2[-\frac b2, \frac b2]$. It can be shown that  ${\mathcal H}\{\phi_k(x)\}= i\sgn(k)  \phi_k(x)$, which implies that functions $\{\phi_k(x)\}$ are the eigenfunctions of the Hilbert transform ${\mathcal H}$ with associated eigenvalues $\pm i$ with $x\in[-\frac b2,\frac b2]$. For a causal periodic continuation, according to (\ref{Hilbert_pair}), we need $\Im{\mathcal C}(H)(x)$ to be the  Hilbert transform of $-\Re{\mathcal C}(H)(x)$. It can be shown \cite{Barannyk_Aboutaleb_Elshabini_Barlow_IEEE} that this implies  $\alpha_{k}=0$ for $k\leq 0$ in  (\ref{E1}). Hence, a causal Fourier continuation has the form
\begin{equation} \label{E3_0}
{\mathcal C}(H)(x)=\sum_{k=1}^{M}  \alpha_k \phi_k(x). 
\end{equation}
Evaluating $H(x)$ at points $x_j$, $j=1,\ldots,\tilde N$, $x_j\in[-0.5, 0.5]$, produces a complex valued system
\begin{equation} \label{E3}
\sum_{k=1}^{M}  \alpha_k \phi_k(x_j)= H(x_j)
\end{equation}
with $\tilde N$ equations for $M$ unknowns $\alpha_{k}$, $k=1,\ldots,M$, $\tilde N\geq M$. If $\tilde N>M$,  the system (\ref{E3})  is overdetermined and has to be solved in the least squares sense. When Fourier coefficients $\alpha_k$ are computed, formula (\ref{E3_0}) provides reconstruction of $H(x)$ on $[-0.5, 0.5]$. The least squares problem is extremely ill-conditioned. However, it can be regularized using a truncated SVD method when singular values below some cutoff tolerance $\xi$  close to the machine precision are being discarded. 
To have a better control on ill-conditioning of matrix problem (\ref{E3}),   more data points $\tilde N$ than the Fourier coefficients $M$ should be used.
%
%
We use at least $\tilde N=2M$ as an effective way to obtain an 
accurate and reliable approximation of $H(x)$ over the  interval $[-0.5, 0.5]$. This relation corresponds\footnote{In \cite{Barannyk_Aboutaleb_Elshabini_Barlow_IEEE}, $N$ denoted the number of points on $[-0.5, 0.5]$, while in this work $N$ is the number of points on $[0, 0.5]$ or originally on $[w_{min},w_{max}]$.} to $N=M$, where $N$ is the number of data points available originally on $[w_{min},w_{max}]$.

Since $\Re H(x)$ and $\Im H(x)$ are even and odd functions of $x$, respectively, the Fourier coefficients 
\[
\alpha_k=\frac{1}{b}\int_{-b/2}^{b/2} H(x) \overline{\phi_k(x)}dx, \quad k=1,\ldots, M,
\]
are real. Here 
$\bar{ \  }$ denotes the complex conjugate. 
To ensure that  numerically computed Fourier coefficients $\alpha_k$ are real, instead of solving complex-valued system (\ref{E3}), one  can separate the real and imaginary parts of ${\mathcal C}(H)(x_j)$ to obtain real-valued system
\begin{equation} \label{E4}
\begin{array}{l}
\displaystyle
\phantom{-}
\Re{\mathcal C}(H)(x_j)=\sum_{k=1}^{M}  \alpha_{k}   \Re\phi_k(x_j), \\[13pt]
\displaystyle
\hspace{10pt}
\Im{\mathcal C}(H)(x_j)=\sum_{k=1}^{M}  \alpha_{k} \Im\phi_k(x_j).
\end{array}
\end{equation}
We show in \cite{Barannyk_Aboutaleb_Elshabini_Barlow_IEEE} that real formulation (\ref{E4}) provides slightly more accurate results than complex.


To access the quality of approximation of $H(x)$ with its causal Fourier continuation ${\mathcal C}(H)(x)$, we introduce reconstruction errors $E_R(x)$ and $E_I(x)$, 
\begin{equation} \label{error_Re}
E_R(x)=\Re H(x) - \Re{\mathcal C}(H)(x),
\end{equation}
\begin{equation}\label{error_Im}
E_I(x)=\Im H(x) - \Im{\mathcal C}(H)(x)
\end{equation}
on the original interval $[-0.5,0.5]$. 

The error analysis performed in \cite{Barannyk_Aboutaleb_Elshabini_Barlow_IEEE} (see also Appendix) shows that the error between $H(x)$ and its causal Fourier continuation $ {\mathcal C}(H+\varepsilon) $ under the presence of a noise $\varepsilon$, has the following upper bound:
\begin{equation} \label{error_estimate}
||H- {\mathcal C}(H+\varepsilon) ||_{L_2(\Omega)} \leq \epsilon_F+\epsilon_n + \epsilon_T.
\end{equation}
Here
\begin{equation}\label{error_Fourier}
\epsilon_F= (1+\Lambda_2 \sqrt{2N(M-K)}) || H-\hat H_{M} || _{L_\infty(\Omega)} 
\end{equation}
is the error due to approximation of $H$ with a causal Fourier series and it decays as ${\mathcal O}(M^{-k+1})$, where $k$ is the smoothness order of the transfer function $H(x)$.
\begin{equation}\label{error_truncation}
\epsilon_T=\Lambda_1 \sqrt{K/b} ||\hat H_{M}||_{L_\infty(\Omega^c)}
\end{equation} 
is the error due to the truncation of singular values and it is typically small and close to the cut-off value $\xi$. 
As (\ref{error_truncation}) indicates, $\epsilon_T$ depends on $b$ and the function $H$ being approximated. 
\begin{equation}\label{error_noise}
\epsilon_n=(1+\Lambda_2 \sqrt{2N(M-K)}) ||\varepsilon||_{L_\infty(\Omega)} 
\end{equation} 
is the error due to the presence of a noise or approximation errors in the given data and it shows a level of causality violation. In practice the size of $\epsilon_n$ is close to the size of noise in data. Function $\hat H_{M}$ and  constants  $\Lambda_1$,  $\Lambda_2$  and $K$ are defined in Appendix. These constants depend only on the continuation parameters $N$, $M$, $b$ and $\xi$ as well as location of discrete points $x_j$, and not on the function $H$. 

The error bound (\ref{error_estimate}) shows that the reconstruction errors $E_R$ and $E_I$ decrease as $M$ increases due to the causal Fourier series approximation error with the error bound $\epsilon_F$ 
 until either the level $\epsilon$ of a noise or level $\epsilon_T$ due to truncation of singular values is reached. If only round-off errors are present in data, the errors will level off at $\epsilon_T$. If reconstruction errors level off at some value $\epsilon>\epsilon_T$ as the resolution increases, the data are declared non-causal with the error approximately at the order of $\epsilon$. More information about the error analysis for the causality characterization methods based on causal Fourier continuations can be found in \cite{Barannyk_Aboutaleb_Elshabini_Barlow_IEEE}. 

%
 
\section{Time Delay Estimation} \label{delay_estimation}

The above approach for causality assessment can be transformed into a delay estimation algorithm by observing the following. Suppose that $h(t)$ is non-zero only from time $T_0\geq 0$, and we would like to identify the time delay $T_0$. Consider the Fourier transform $H(w)$ of $h(t)$:
\[
H(w)=\int_{t=T_0}^\infty \e^{-iwt} h(t) dt=\int_{t=T_0}^\infty \e^{-i\frac x a t} h(t) dt
\]
where we used the substitution $x=aw$, $a=\frac{0.5}{w_{max}}$. Introducing $\tau=\frac ta$, we can write
\[
H(w)=a \int_{\frac{T_0}{a}}^\infty \e^{-ix\tau} h(a\tau) d\tau,
\]
or with $u=\tau-\frac{T_0}{a}$, we obtain
\[
H(w)=a \e^{-ix\frac{T_0}{a}} \int_0^\infty \e^{-ixt} h(T_0+au) du = a \e^{-ix\frac{T_0}{a}}  G(x),
\]
where $G(x)$ is the Fourier transform of a causal function with no time delay.  This implies that when $0\leq T \leq {T_0}$, the transfer function $H(x)\e^{ix\frac Ta}$ is causal, but when $T\geq {T_0}$, the transfer function $H(x)\e^{ix\frac Ta}$ has a non-causal component. Therefore, $\tilde T_0={T_0}/a$ is the time delay for $H(x)$, and the delay $T_0$ for the original function $H(w)$ is recovered by multiplying $\tilde T_0$ by $a$. Since one can add any integral multiple of $2\pi$ to $xT/a$, it is enough to restrict our investigations to the interval
\[
0\leq \frac Ta \leq T_{max}=\frac{2\pi}{x_{max}}=\frac{2\pi}{0.5}=4\pi.
\]
Then for each potential time delay $0\leq \frac Ta\leq T_{max}$, we solve the following modified system
\begin{equation} \label{ED3}
\sum_{k=1}^{M}  \alpha_k \phi_k(x_j)=\e^{ix_j\frac Ta}H(x_j), \quad j=1,\ldots,\tilde N
\end{equation}
or its equivalent real-valued formulation.
For $T< T_0$, the reconstruction errors $E_R$, $E_I$ should be small and approximately of the same order. As $T$ increases and becomes greater than some critical transition time close to the time delay $ T_0$, the reconstruction errors should start increase. 
The goal is to approximate $T_0$. The difficulty is that the reconstruction errors grow gradually as $T\geq  T_0$, so transition is not sharp. Moreover, the order of reconstruction errors for $T<  T_0$ depends on the resolution of data and threshold $\xi$ used in the truncated SVD method, which, in turn, affects a transition time. In addition, 
a noise in data, if present, also affects when reconstruction errors start growing. 
A similar approach was used in \cite{Knockaert_Dhaene_2008} to estimate the time delay for square integrable transfer functions. In this contribution, we extend the approach to more general transfer functions. In addition, we use a different causality measure than in \cite{Knockaert_Dhaene_2008} and take into account different resolutions of given data and a possible presence of noise.
The approach can be extended to multi-port and mixed mode networks by applying it to each element of the transfer matrix.


\section{Numerical Examples} \label{numerical_examples}

In this section, we apply a proposed technique to several analytic and simulated examples when the time delay is either known exactly or can be estimated using other techniques. We also consider the effect of noise presence on the accuracy of timed delay estimation.

\subsection{Four-Pole Example} \label{Example_four_pole}

\begin{figure}[h] \begin{center}
\includegraphics[width=2.5in,angle=0]{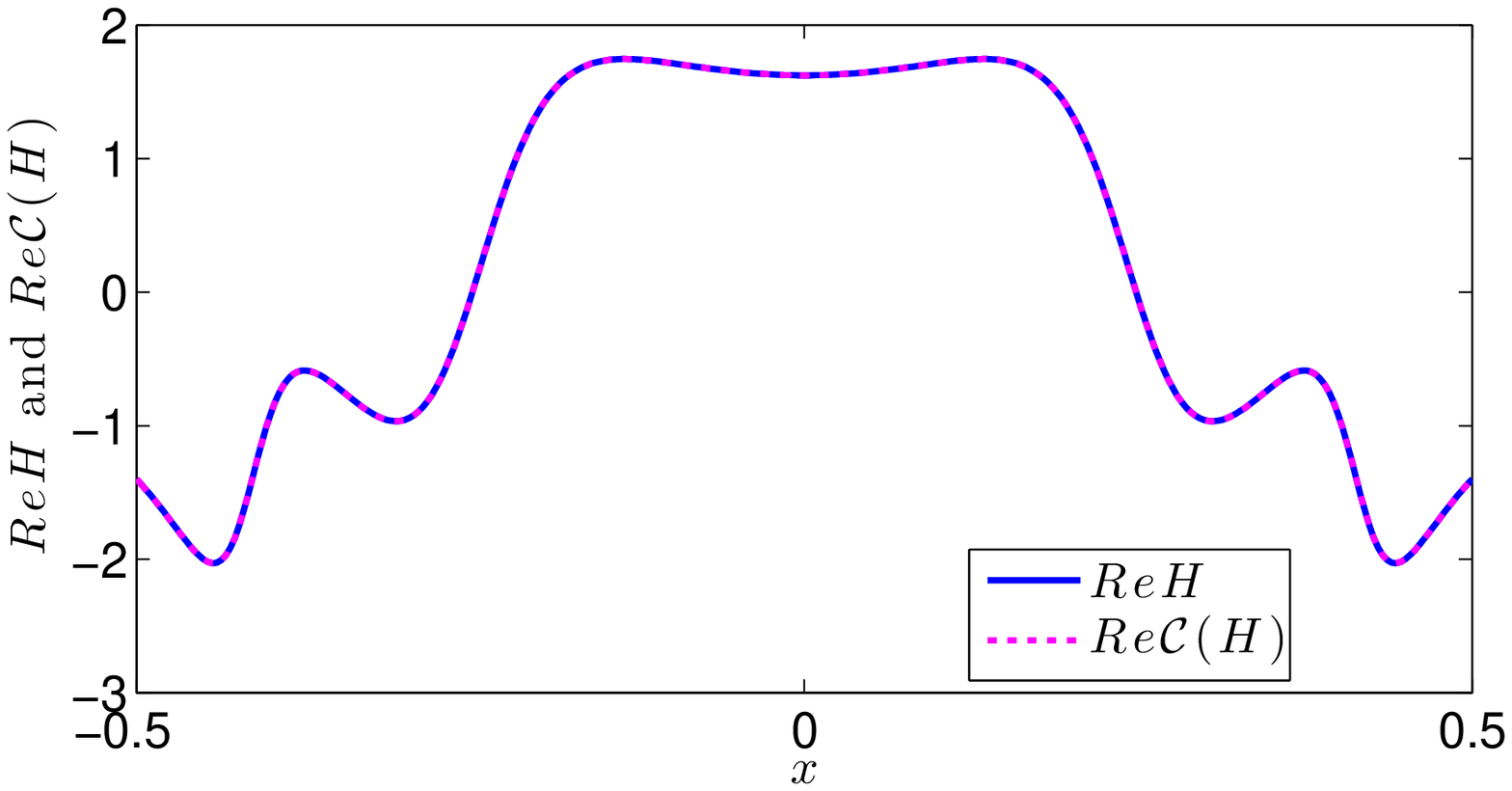}
\hspace{10pt}
\includegraphics[width=2.5in,angle=0]{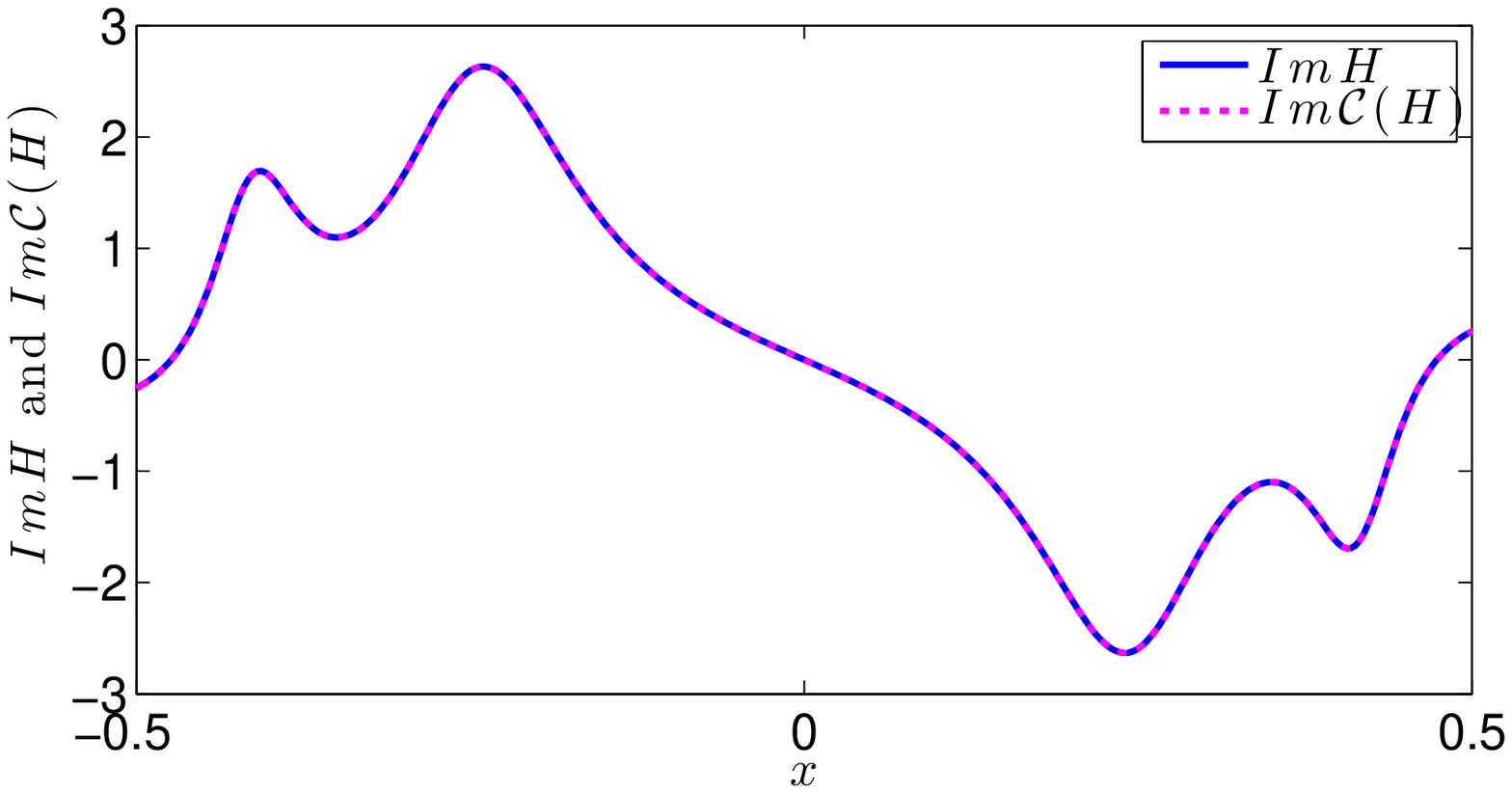}
\end{center}
\caption{$\Re H(x)$ and $\Im H(x)$ of the rescaled and reflected by symmetry  transfer function $H(x)$ in the four-pole example together with their causal Fourier continuations $\Re {\mathcal C}(H)$ and $\Im {\mathcal C}(H)$, respectively, with $M=300$ Fourier coefficients.} 
\label{Ffour_pole0_2}
\end{figure}

Consider a transfer function with four poles and time delay $T_0$, defined by
\begin{equation} \label{Exfour_pole1}
H(w)=\e^{-iwT_0} \tilde H(w)
\end{equation}
with
\[
\tilde H(w)=\frac{r_1}{iw+p_1}+\frac{\bar r_1}{iw+\bar p_1}+\frac{r_2}{iw+p_2}+\frac{\bar r_2}{iw+\bar p_2} 
\]
%
where $r_1=1+2i$, $p_1=1+3i$, $r_2=\frac 23+\frac 12 i$, $p_2=\frac 12+5i$, and  $T_0=0.25$. Since the poles of $ \tilde H(w)$  are located in the upper half $w$-plane at $\pm 3+i$  and $\pm 5+\frac 12 i$ , this function is causal as a sum of four causal transforms, and has no time delay. Therefore, the function $H(w)$ is a  causal function delayed with offset $T_0$.
$H$ is sampled on $[0,w_{max}]$ at $N$ frequency points varying from $50$ to $1500$ with $w_{max}=6$.
\begin{figure}[h] \begin{center}
\includegraphics[width=2.5in,angle=0]{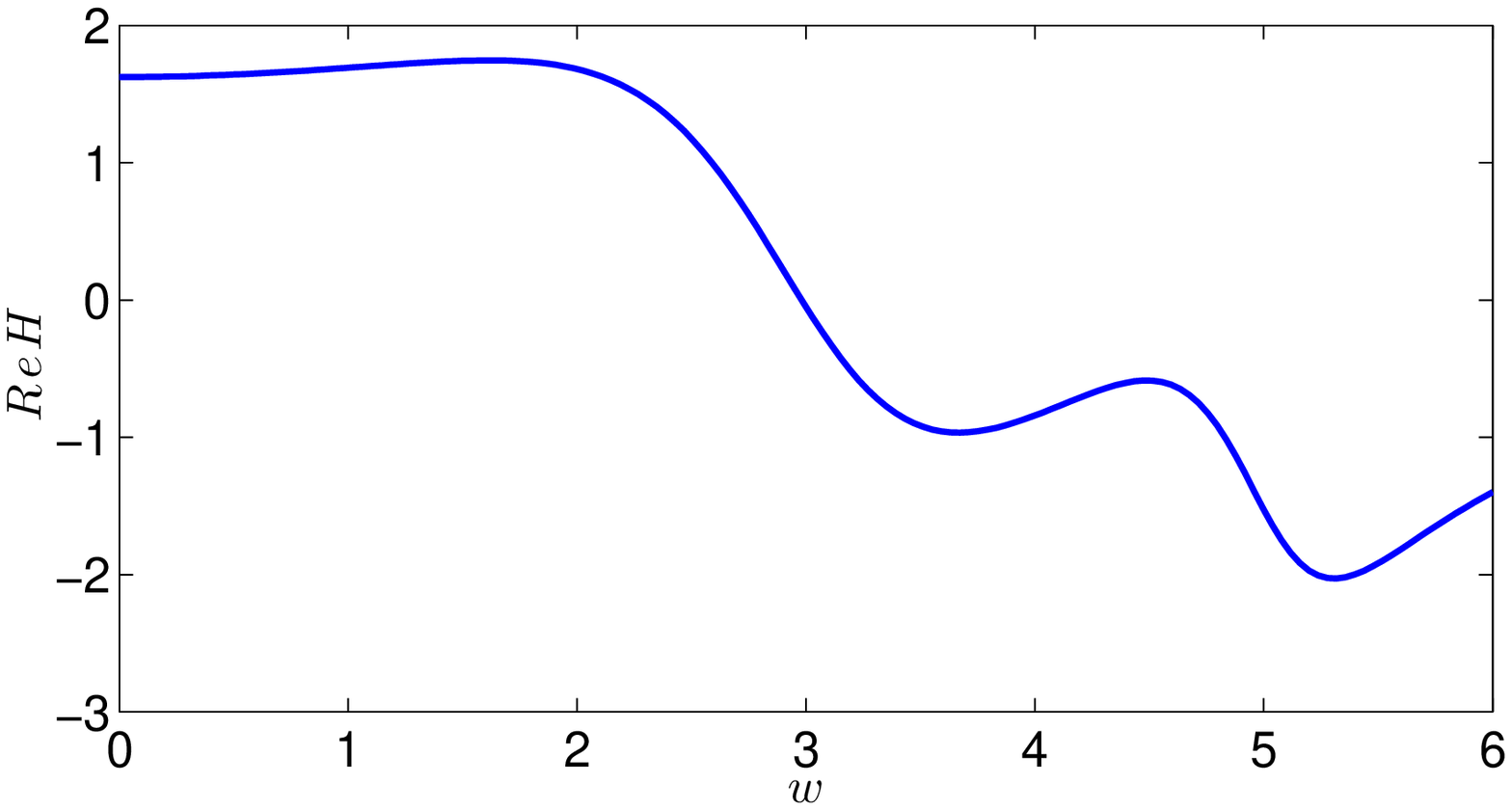}
\hspace{10pt}
\includegraphics[width=2.5in,angle=0]{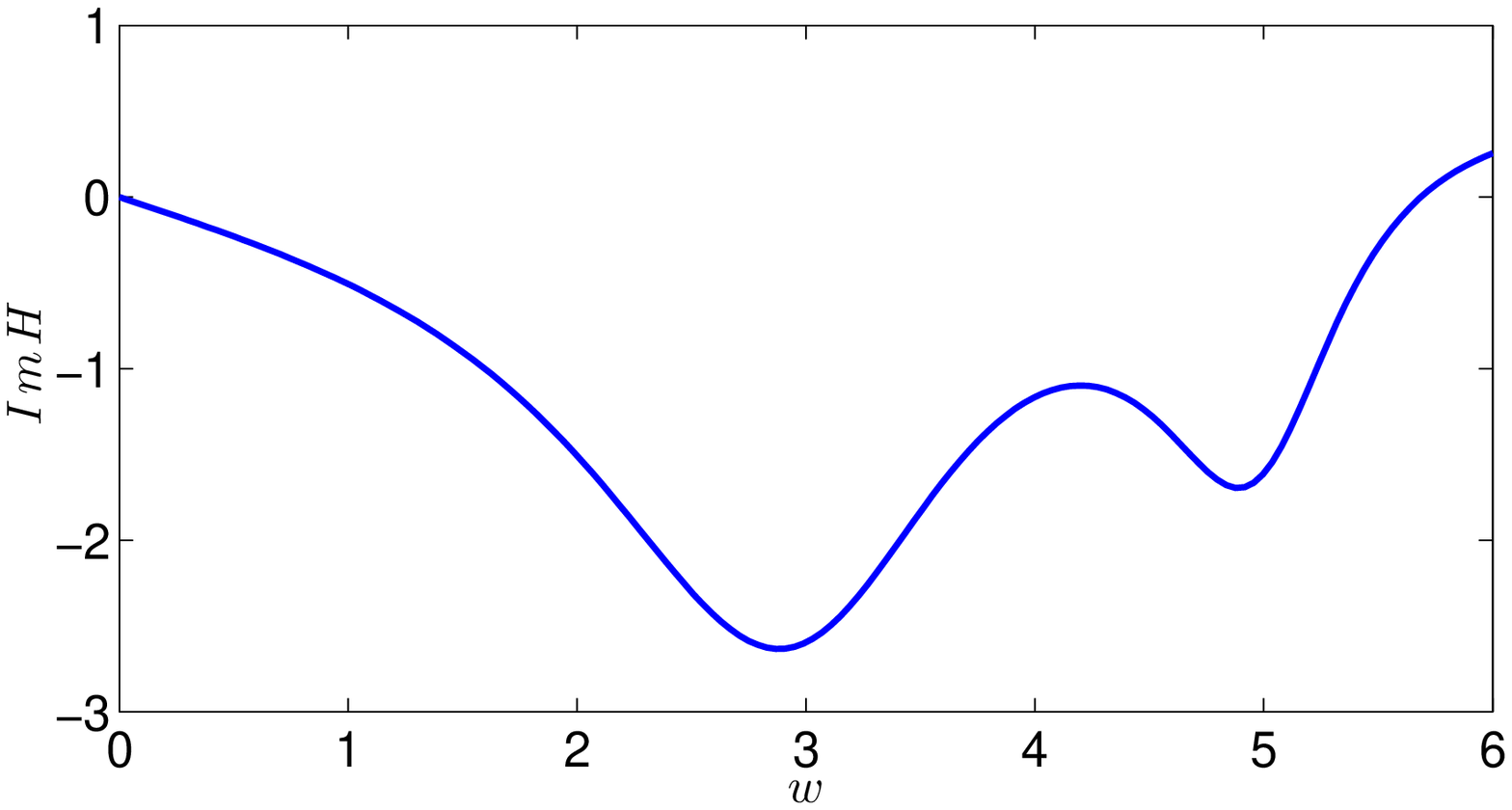}
\end{center}
\caption{$\Re H(w)$ and $\Im H(w)$ in the four-pole example.} 
\label{Ffour_pole0}
\end{figure}
The real and imaginary parts of $H(w)$ are shown in Fig. \ref{Ffour_pole0}. After rescaling with  $x=\frac{0.5}{w_{max}}w$ and reflecting to negative frequencies, we obtain a rescaled  transfer function $H(x)$  defined on  $x\in[-0.5, 0.5]$, for which we construct a causal Fourier continuation ${\mathcal C}(H)$ defined in (\ref{E3_0}) using $M=N$ Fourier coefficients. 
Hence, the number $M$ of Fourier coefficients also varies between $50$ and $1500$. $\Re H(x)$    and $\Im H(x)$ of the rescaled and reflected $H(x)$ together with their causal Fourier continuations with $M=300$ are  depicted in Fig. \ref{Ffour_pole0_2}. Even though given $H(x)$ and its causal Fourier continuation approximation look indistinguishable, the actual reconstruction errors $E_R$ and $E_I$ in both real and imaginary parts, that are defined in (\ref{error_Re}), (\ref{error_Im}),   are on the order of $10^{-6}$ and they decreases as $M$ increases (with $M=N$). For example, with $M=800$, the errors are on the order of $10^{-13}$. Since both errors $E_R$ and $E_I$ are of the same order, it is enough to analyze one of the errors, for example,  $E_R$. The results using  $E_I$ are similar.

\begin{figure}[h] \begin{center}
\includegraphics[width=2.8in,angle=0]{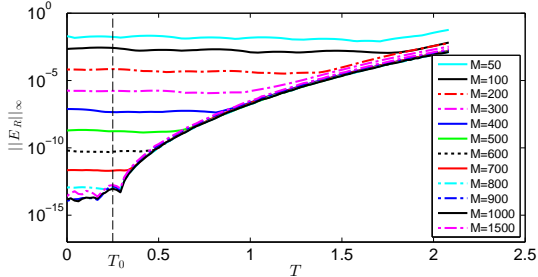}
\end{center}
\caption{Evolution of the reconstruction error $||E_R||_\infty$ as a function of $T$ in the four-pole example. The dashed line corresponds to the exact delay $T_0=0.25$.} 
\label{Ffour_pole1}
\end{figure}

\begin{figure}[h] \begin{center}
\includegraphics[width=2.5in,angle=0]{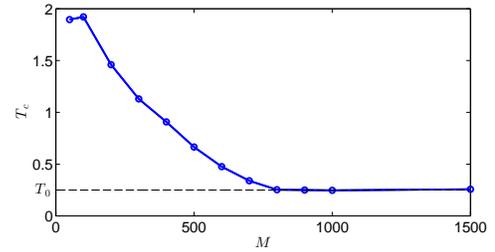}
\end{center}
\caption{Critical transition times $T_c$ in the  four-pole example that approach  $T_0$ as $M$ increases. The dashed line corresponds to the exact delay $T_0=0.25$.} 
\label{Ffour_pole2}
\end{figure}

To estimate the time delay, we analyze the evolution of the $||E_R||\infty$, shown in Fig. \ref{Ffour_pole1}, for various values $M$. 
Since the error due to a causal Fourier series approximation decreases with $M$ (see error bound (\ref{error_Fourier})), the reconstruction error between the given transfer function $H$ and its causal Fourier continuation ${\mathcal C}(H)$  also  decreases as $M$ increases until it  either reaches the level $\xi$ of filtering of singular values or a level $\epsilon$ of noise/causality violations (see error bounds (\ref{error_truncation}) and  (\ref{error_noise}), respectively). 
%
%
%
%
%
%
For each fixed $M$, as time $T$ increases, the errors $E_R$ and $E_I$ first are small and about of the same order until some transition time close to the time delay $T_0$ is approached. After that the errors  grow approximately as a power function on the loglog scale. For smaller $T$, the errors are dominated by a causal Fourier series approximation error and then for $T$ greater than some transition time -- by causality violations since this value $T$ provides large enough negative time delay and shifts a causal function into a non-causal area.
A transition value $T=T_c$, we call it a critical  time, from a plateau region to a growth region, is different for each $M$ and it decreases as the resolution or number of Fourier coefficients increases if the error is dominated by the causal Fourier series approximation error. The critical times $T_c$ approach the time delay $T_0$ as $M$ increases. The goal is to estimate $T_0$ using the error curves  shown in Fig. \ref{Ffour_pole1}.
Analyzing graphs of the error curves  for $M>800$, we observe some non-monotonic behavior at $T$ close to $T_0$. This behavior is due to the filtering of the singular values below the threshold $\xi=10^{-13}$ that we used in our experiments. By increasing the value of $\xi$, the non-monotonic behavior will be present at smaller values of $M$. This suggests that portions of error curves close to threshold $\xi$ are affected by filtering and may be inaccurate and difficult to use for time delay estimation as we  find in our experiments.
%
%
To estimate critical times $T_c$ of transition from the plateau region to the growth region, we approximate the growing region by  a quadratic function on the loglog scale.
Specifically, we assume that 
$\ln T\approx a_2 \left(\ln ||E_R||_\infty\right)^2+a_1 \ln ||E_R||_\infty +a_0\equiv f(\ln ||E_R||_\infty)$, where coefficients $a_0$, $a_1$, and $a_2$  are determined in the least squares sense. The resulting quadratic function $f(\ln ||E_R||_\infty)$ is then evaluated at the value of $||E_R||_\infty$ at $T=0$ that is assumed to be the ``most causal" time.  By taking exponential function of the result, we find a critical transition time $T_c$ for a given $M$. This procedure produces estimates of the time delay $T_0$ for various values of $M$. The graph of the critical transition times $T_c$ as a function of $M$ is shown in Fig. \ref{Ffour_pole2}. One can clearly see that the critical times approach the exact time delay $T_0=0.25$ as $M$ increases. A good approximation of $T_0$ is achieved at $M=800$. 
%
The values of $T_c$ for $M\geq 200$ are presented in Table \ref{Tcrit_four_pole}. The results indicate that the approximations become more accurate as $M$ increases. The error with $M=900$  is less than 1\,\%. At the same time, the error with $M=1500$ is about 3\%, which is due to the fact that the results in this case are more affected by the filtering of singular values. In the cases when $M$ is high and the resulting error is not flat for $T<T_0$, instead of evaluating a fitted quadratic curve at the value of $||E_R||_\infty$ at $T=0$ we evaluate it at $\xi$, the threshold of filtering singular values, to avoid using results affected by filtering.
\begin{table}[h]
\begin{center}
\begin{tabular}{|c|c||c|c|}
\hline
\rule{0cm}{10pt}
$M$ & $T_c$ & $M$ & $T_c$  \\[3pt]
\hline
\rule{0cm}{10pt}
200 &    1.4604 & 700  & 0.3394  \\[3pt]
\hline
\rule{0cm}{10pt}  
300 &  1.1294  &   800  &  0.2529 \\[3pt]
\hline
\rule{0cm}{10pt}  
400  &  0.9077 & 900  & 0.2497 \\[3pt] 
\hline
\rule{0cm}{10pt}    
500  &  0.6655 & 1000  &  0.2472 \\
\hline
\rule{0cm}{10pt}    
600  &  0.4759 & 1500  &  0.2576 \\
\hline
\end{tabular}
\end{center}
\caption{Critical transition times $T_c$ in the four-pole example that approach $T_0=0.25$  as $M$ increases.}
\label{Tcrit_four_pole}
\end{table} 
       
\begin{figure}[h] \begin{center}
\includegraphics[width=2.7in,angle=0]{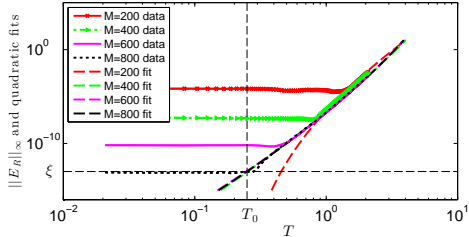}
\end{center}
\caption{$||E_R||_\infty$ in the four-pole example with $M=200$, 400, 600 and 800  together with their extrapolated quadratic fits. Vertical dashed line indicates the exact time delay $T_0=0.25$, while horizontal dashed line indicates the level of filtering of singular values given by $\xi=10^{-13}$.}
\label{Ffour_pole3}
\end{figure}

In practice the number $N$ of samples of the transfer function $H(w)$ is usually limited, which sets the bound for number $M=N$ of Fourier coefficients, so it may not always be possible to use large enough $M$ to obtain critical time $T_c$ close enough to the actual time delay $T_0$. A good method should be capable of producing an accurate approximation of $T_0$ even with a small number of data points. We achieve this by employing another approach for time delay estimation. Using the obtained   fitted quadratic error curves, we extrapolate them  to the value $\xi$ of filtering of singular values, which is typically chosen to be close to the machine precision. This corresponds to finding time $T$ at which the error reaches the value $\xi$. This choice is natural since the errors below $\xi$ are most likely affected by filtering and may not be accurate enough to use. The results of such extrapolation are shown in Fig. \ref{Ffour_pole3} for   $M=200$, 400, 600 and 800. An intersection of the extrapolated curve corresponding to $M=200$ is at a value $T=0.45451$, which is a bit far from the exact $T_0=0.25$. At the same time, intersections of extrapolated curves with higher values of $M$ are  much closer -- see  Table \ref{T0_extrapol_four_pole} for details. 
\begin{table}[h]
\begin{center}
\begin{tabular}{|c|c||c|c|}
\hline
\rule{0cm}{10pt}
$M$ & $T_0$ approximation & $M$ & $T_0$ approximation  \\[3pt]
\hline
\rule{0cm}{10pt}
200 &    0.45451 & 700  & 0.24734  \\[3pt]
\hline
\rule{0cm}{10pt}  
300 &  0.27235  &   800  &  0.24969 \\[3pt]
\hline
\rule{0cm}{10pt}  
400  &  0.25297 & 900  & 0.24974 \\[3pt] 
\hline
\rule{0cm}{10pt}    
500  &  0.25053 & 1000  &  0.24724 \\
\hline
\rule{0cm}{10pt}    
600  &  0.24633 & 1500  &  0.25759 \\
\hline
\end{tabular}
\end{center}
\caption{Approximations of $T_0=0.25$ in the four-pole example using extrapolation.}
\label{T0_extrapol_four_pole}
\end{table}
%
Results shown in this  table 
indicate that as $M$ increases, extrapolated quadratic curve intersect the horizontal line the value $\xi$ at times closer to $T_0$.  Obtained approximations of $T_0$ 
can be averaged producing $T_0\approx 0.24805$. The approach with extrapolation provides a faster convergence and good approximations of $T_0$ even for small values of $M$, i.e. less data points are needed to approximate $T_0$.

\begin{figure}[h] \begin{center}
\includegraphics[width=2.5in,angle=0]{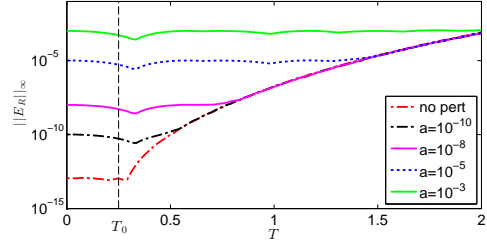}
\end{center}
\caption{Evolution of $||E_R||_\infty$ in the four-pole example with added sine perturbation $a\sin(10\pi x)$. The dashed line corresponds to the exact delay $T_0=0.25$.} 
\label{Ffour_pole4}
\end{figure}

\begin{figure}[h] \begin{center}
\includegraphics[width=3.0in,angle=0]{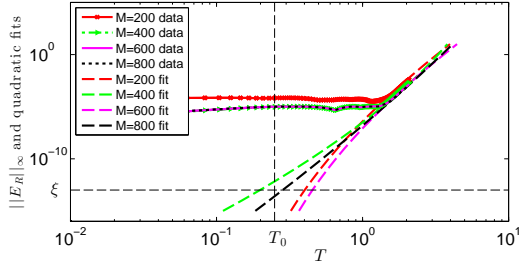}
\end{center}
\caption{Error curves $||E_R||_\infty$ in the four-pole example  with  the added perturbation $10^{-5}\sin(10\pi x)$  and $M=200$, 400, 600 and 800  together with their extrapolated quadratic fits. Vertical dashed line indicates the exact time delay $T_0=0.25$, while horizontal dashed line indicates the level of filtering of singular values given by $\xi=10^{-13}$.}
\label{Ffour_pole6}
\end{figure}

\begin{figure}[h] \begin{center}
\includegraphics[width=3.0in,angle=0]{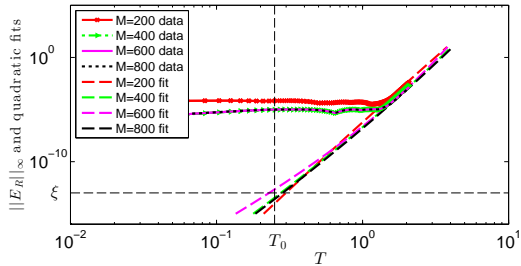}
\end{center}
\caption{$||E_R||_\infty$ in the four-pole example  with  the added perturbation $10^{-5}\sin(10\pi x)$  and $M=200$, 400, 600 and 800  together with their extrapolated quadratic fits constructed using more narrow fitting region. Vertical dashed line indicates the exact time delay $T_0=0.25$, while horizontal dashed line indicates the level of filtering of singular values given by $\xi=10^{-13}$.}
\label{Ffour_pole6_2}
\end{figure}

We also consider the effect of noise on the time delay estimation. To study this, we impose a sine perturbation
\begin{equation} \label{sine_pert}
a\sin(10\pi x)
\end{equation}
of various amplitudes $a$ that we add to $\Re H$, while keeping $\Im H$ unchanged. We choose $N=800$ 
and vary $a$ from $10^{-10}$ to $10^{-3}$. The reconstruction error $E_R$ with no perturbation for early times $T<T_0$ is of the order of $10^{-13}$, as shown in Fig. \ref{Ffour_pole4}, that corresponds to the level of filtering of singular values. When the perturbation is added, the reconstruction errors for $T<T_0$ are higher and approximately of the order of $a$. Once some critical transition time greater than $T_0$ is passed, reconstruction errors start growing and they grow at the same rate and coincide almost perfectly with each other. This observation suggests that the proposed approach can also be  used  in the cases when data have a noise, which is typical in real-life applications, when data have either measurement or simulation errors. 
For noise with a smaller amplitude, the region close to $T_0$ will be less affected by noise and a bigger growing region will be available for fitting, so we expect better accuracy of time delay estimation in such cases. 
When more noise in data is present, less growing region will be available for fitting and extrapolation of fitted quadratic error curves may be less accurate.  We demonstrate this by considering two cases: with $a=10^{-5}$ (noisier case) and $a=10^{-8}$ (less noisy case).
\begin{figure}[h] \begin{center}
\includegraphics[width=2.8in,angle=0]{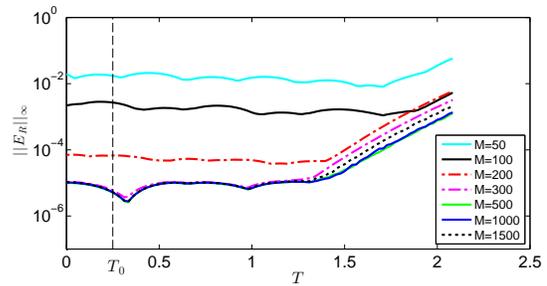}
\end{center}
\caption{Evolution of  $||E_R||_\infty$ in the four-pole example with the added  perturbation $10^{-5}\sin(10\pi x)$. The dashed line corresponds to the exact delay $T_0=0.25$.} 
\label{Ffour_pole5}
\end{figure}
The error curves with a higher amplitude $a=10^{-5}$ are presented in Fig. \ref{Ffour_pole5}. It is clear that the error does not become smaller than $10^{-5}$ as $M\geq 300$ gets larger because of the noise. We use available  growing regions and extrapolate fitted error curves to find their intersection with the horizontal line with value $\xi$. This gives us time $T$ when the error reaches the value $\xi$ for each considered $M$. The results of such extrapolation for $M=200$, 400, 600 and 800 are shown in Fig. \ref{Ffour_pole6}.
Clearly, extrapolated error curves reach value $\xi$ at times around $T_0$ but not close enough to $T_0$ and without established convergence but rather in a spread-out manner around $T_0$.
\begin{table}[h]
\begin{center}
\begin{tabular}{|c|c||c|c|}
\hline
\rule{0cm}{10pt}
$M$ & $T_0$ estimate & $M$ & $T_0$ estimate \\[3pt]
\hline
\rule{0cm}{10pt}
  200 &    0.40158 & 700  & 0.39113 \\[3pt]
\hline
\rule{0cm}{10pt}  
  300 &  0.25578  &   800  &  0.28392 \\[3pt]
\hline
\rule{0cm}{10pt}  
400  &  0.19863 & 900  & 0.26543 \\[3pt] 
\hline
\rule{0cm}{10pt}    
500  &  0.14311 & 1000  &  0.20293 \\
\hline
\rule{0cm}{10pt}    
600  &  0.45358 & 1500  &  0.32837 \\
\hline
\end{tabular}
\end{center}
\caption{Approximations of $T_0$ in the four-pole example with perturbation  $10^{-5}\sin(10\pi x)$ using extrapolations with original fitting regions for various $M$. The exact value $T_0=0.25$, averaged value $T_0^{(aver)}=0.26586$.}
\label{T0approx_ampn5_table}
\end{table} 
%
Approximations of $T_0$ for values of $M$ that we investigated are shown in Table \ref{T0approx_ampn5_table}.
Averaging these approximations we obtain  $T_0^{(aver)}=0.26586$. The extrapolated curves can be made more focused around $T_0$ by narrowing down the fitted region. Results of this procedure are  shown in Fig. \ref{Ffour_pole6_2}. This improves a little an average time delay to $T_0^{(aver)}=0.24216$.
%

%

\begin{figure}[h] \begin{center}
\includegraphics[width=2.8in,angle=0]{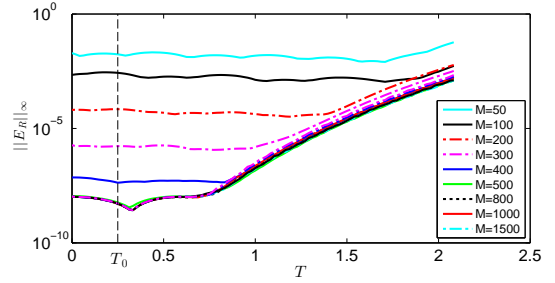}
\end{center}
\caption{Evolution of $||E_R||_\infty$  in the four-pole example with the added  perturbation $10^{-8}\sin(10\pi x)$ of smaller amplitude. The dashed line corresponds to the exact delay $T_0=0.25$.} 
\label{Ffour_pole7}
\end{figure}

\begin{figure}[h] \begin{center}
\includegraphics[width=3.0in,angle=0]{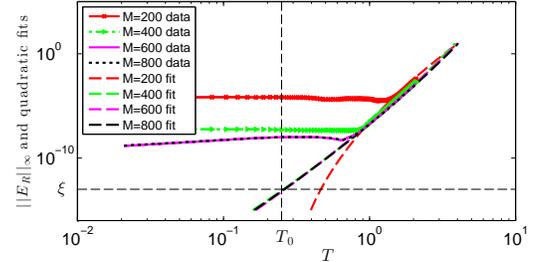}
\end{center}
\caption{$||E_R||_\infty$ in the four-pole example  with  the added perturbation $10^{-8}\sin(10\pi x)$  and $M=200$, 400, 600 and 800  together with their extrapolated quadratic fits. Vertical dashed line indicates the exact time delay $T_0=0.25$, while horizontal dashed line indicates the level of filtering of singular values given by $\xi=10^{-13}$.}
\label{Ffour_pole8}
\end{figure}

Next we show results when a smaller noise of amplitude $a=10^{-8}$ is added. The evolution of $||E_R||_\infty$ as $T$ increases is shown in Fig. \ref{Ffour_pole7}. We can see that the plateau error region in this case is at about $10^{-9}$ level, so the error growth region is bigger than in the previous case, which should make fitting and extrapolation more accurate.
Indeed, extrapolated quadratic curves intersect the horizontal line with value $\xi$ in a more localized region about $T_0$ as shown in Fig. \ref{Ffour_pole8}, while averaging of obtained approximation to $T_0$ produces $T_0^{(aver)}=0.25436$, which is more accurate than  in the case with a higher amplitude $a=10^{-5}$.

\begin{table}[h]
\begin{center}
\begin{tabular}{|c|c||c|c|}
\hline
\rule{0cm}{10pt}
$M$ & $T_0$ estimate & $M$ & $T_0$ estimate \\[3pt]
\hline
\rule{0cm}{10pt}
  200 &    0.46392 & 700  & 0.25502 \\[3pt]
\hline
\rule{0cm}{10pt}  
  300 &  0.27635  &   800  &  0.26081 \\[3pt]
\hline
\rule{0cm}{10pt}  
400  &  0.25683 & 900  & 0.2444 \\[3pt] 
\hline
\rule{0cm}{10pt}    
500  &  0.26798 & 1000  &  0.25582 \\
\hline
\rule{0cm}{10pt}    
600  &  0.26391 & 1500  &  0.25292 \\
\hline
\end{tabular}
\end{center}
\caption{Approximations of $T_0$ in the four-pole example with perturbation  $10^{-8}\sin(10\pi x)$ using extrapolations for various $M$. The exact value $T_0=0.25$, averaged value $T_0^{(aver)}=0.25436$.}
\label{T0approx_ampn8}
\end{table} 
%

\subsection{Transmission Line Example} \label{transmission_line_new}

\begin{figure}[h] \begin{center}
\includegraphics[width=2.5in,angle=0]{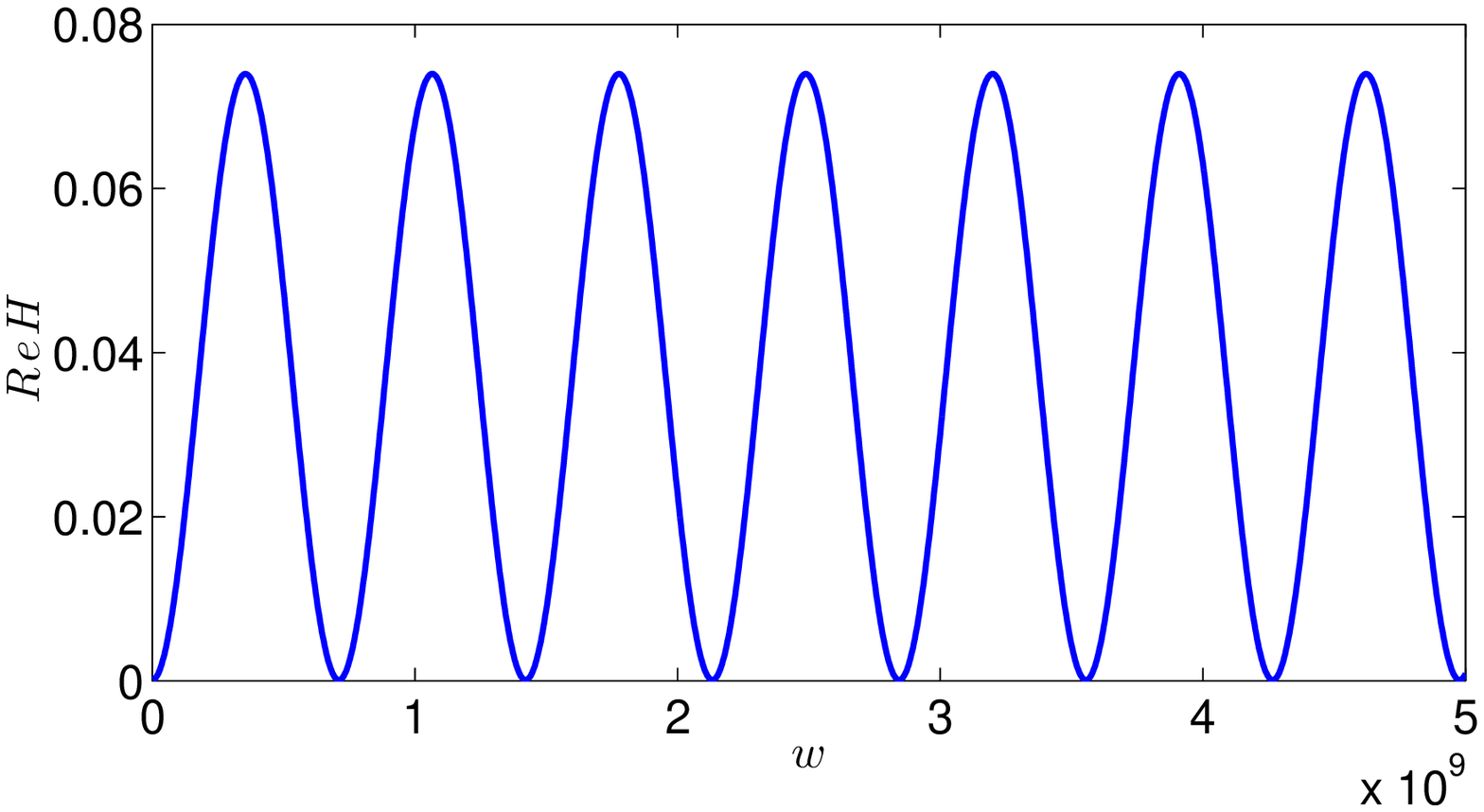}
\hspace{10pt}
\includegraphics[width=2.5in,angle=0]{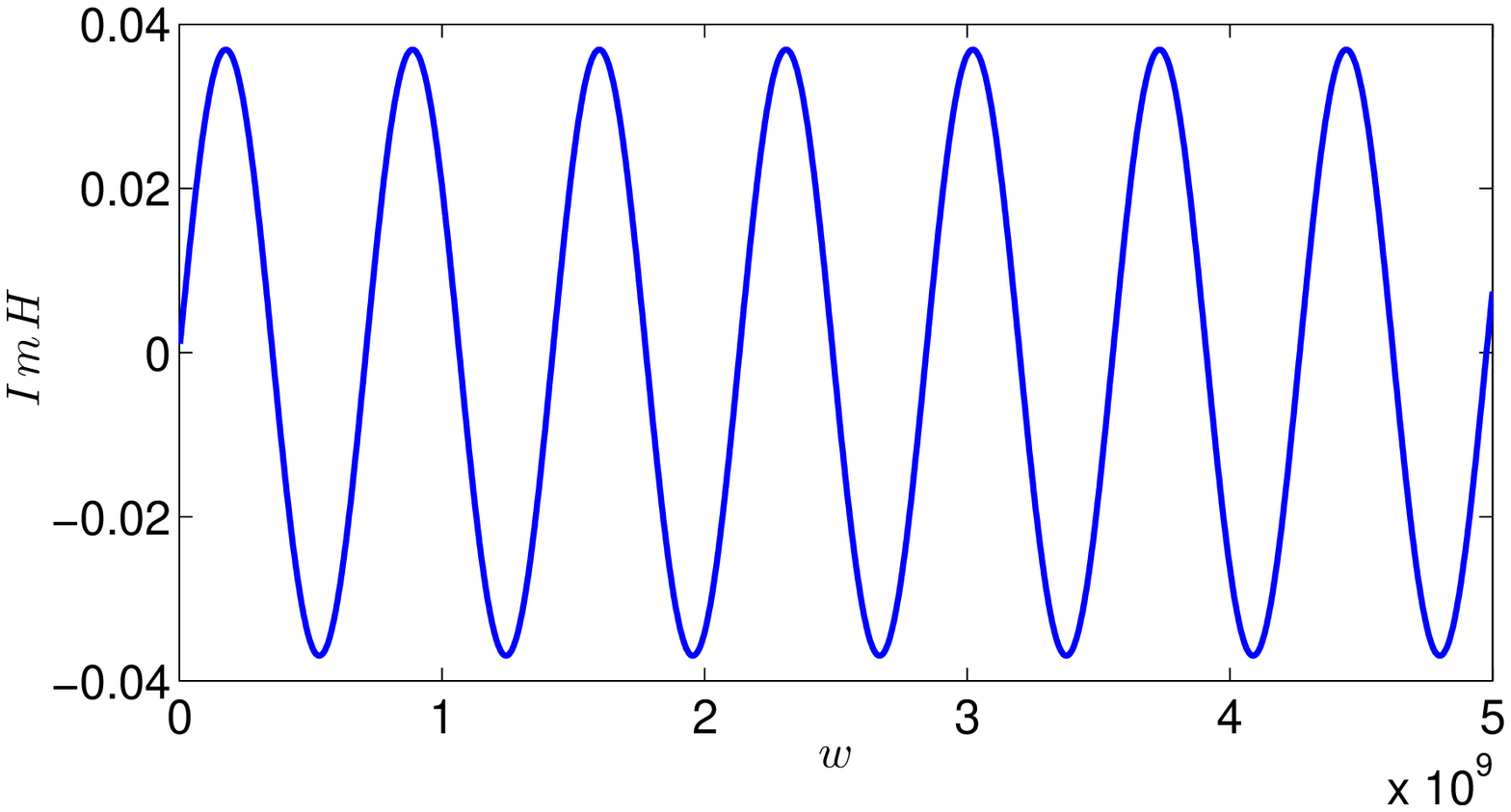}
\end{center}
\caption{$\Re H(w)$ and $\Im H(w)$ in the transmission line example.} 
\label{Ftransm_line_new0}
\end{figure}

\begin{figure}[h] \begin{center}
\includegraphics[width=2.5in,angle=0]{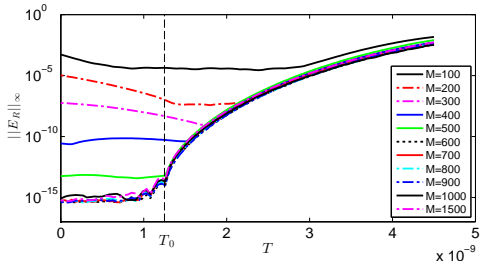}
\end{center}
\caption{Evolution of $||E_R||_\infty$ in the transmission line example as $M$ varies. Vertical dashed line indicates the time delay $T_0=1.25$\,ns.}
\label{Ftransm_line_new1}
\end{figure}

\begin{figure}[h] \begin{center}
\includegraphics[width=2.5in,angle=0]{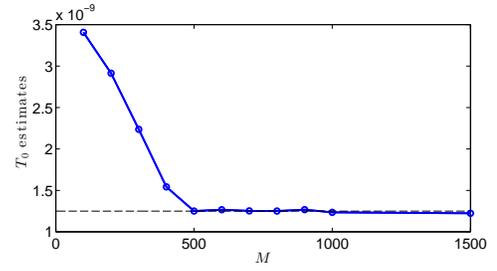}
\end{center}
\caption{Estimation of the delay time in the transmission line example using critical transition times $T_c$ as $M$ varies. The dashed line corresponds to the exact delay $T_0=1.25$\,ns.} 
\label{Ftransm_line_new2}
\end{figure}

\begin{figure}[h] \begin{center}
\includegraphics[width=2.5in,angle=0]{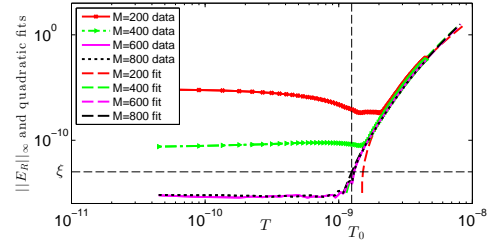}
\end{center}
\caption{$||E_R||_\infty$ in the transmission line example with $M=200$, 400, 600 and 800  together with their quadratic fits. Vertical dashed line indicates the exact time delay $T_0=1.25$\,ns, while horizontal dashed line indicates the level of filtering of singular values given by $\xi=10^{-13}$.} 
\label{Ftransm_line_new3}
\end{figure}

We consider a uniform transmission line segment with the following per-unit-length parameters: $L=7.574$ nH/inch, $C=2.61166$ pF/inch, $R=16.278$ $m\Omega$/inch, $G=5.58$  $\mu$S/inch and length ${\mathcal L}=5$ inches. 
%
%
The frequency is sampled on the interval $(0, 5.0]$ GHz.
The scattering matrix of the structure is computed using Matlab function {\tt rlgc2s}. We consider the element $\tilde H(w)=S_{11}(w)$ and 
impose the time delay $T_0=1.25$\,ns  by multiplying $\tilde H(w)$ by  $\exp(-iwT_0)$ to get the delayed transfer function $H(w)=\exp(-iwT_0)\tilde H(w)$. The real and imaginary parts of $H(w)$ are given in Fig. \ref{Ftransm_line_new0}. 
%
The error curves for different $M$ are shown in Fig. \ref{Ftransm_line_new1} indicating that the reconstruction error decreases quickly with $M$ and reaches the level close to machine precision at $M=600$.
Constructing fitted quadratic error curves and finding their intersections with $||E_R||_\infty$ at $T=0$ or finding times when these fitted error curves reach the value $\xi$ of the error for $M\geq 600$, we get a sequence of critical transition times $T_c$, that we show in Fig. \ref{Ftransm_line_new2}. Clearly, critical times $T_c$ converge to $T_0$ and provide a good approximation of $T_0$ for $M\geq 500$.
Using an alternative approach when we extrapolate the fitted quadratic error curves to find their intersections with the error threshold $\xi$, we find approximations of $T_0$. Some of these curves for $M=200$, 400, 600 and 800 are depicted in Fig. \ref{Ftransm_line_new3}.
Approximations of $T_0$ using extrapolation procedure for various values of $M$ ranging from $M=200$ to $1500$ are given in Table \ref{Ftransm_line_new_table}. A good approximation of $T_0$ is obtained even with $M=300$. As before, approximations of $T_0$ become better as $M$ increases, but for very large values of $M\geq 1000$ when the reconstruction error falls below the filtering threshold $\xi$ and filtering affects the results more, extrapolation becomes less accurate. Averaging obtained approximations of $T_0$ produces $T_0^{(aver)}=1.2498$\,ns, that is very close to the exact value $T_0=1.25$\,ns.
\begin{table}[h]
\begin{center}
\begin{tabular}{|c|c||c|c|}
\hline
\rule{0cm}{10pt}
$M$ & $T_0$ estimate (in ns) & $M$ & $T_0$ estimate (in ns) \\[3pt]
\hline
\rule{0cm}{10pt}
  200 &     $1.5194$ & 700  & $1.2531$ \\[3pt]
\hline
\rule{0cm}{10pt}  
  300 &  $1.3147$  &   800  &  $1.2512$ \\[3pt]
\hline
\rule{0cm}{10pt}  
400  &  $1.2793$ & 900  & $1.2678$ \\[3pt] 
\hline
\rule{0cm}{10pt}    
500  &  $1.2506$ & 1000  &  $1.2348$ \\
\hline
\rule{0cm}{10pt}    
600  &  $1.2668$ & 1500  &  $1.2242$ \\
\hline
\end{tabular}
\end{center}
\caption{Approximations of $T_0$ in the transmission line example using extrapolations for various $M$. The exact value $T_0=1.25$\,ns, averaged value $T_0^{(aver)}=1.2498$\,ns.}
\label{Ftransm_line_new_table}
\end{table} 
%

%

\subsection{Dawson's Integral Example} \label{Dawson_integral_example}

\begin{figure} \begin{center}
\includegraphics[width=2.7in,angle=0]{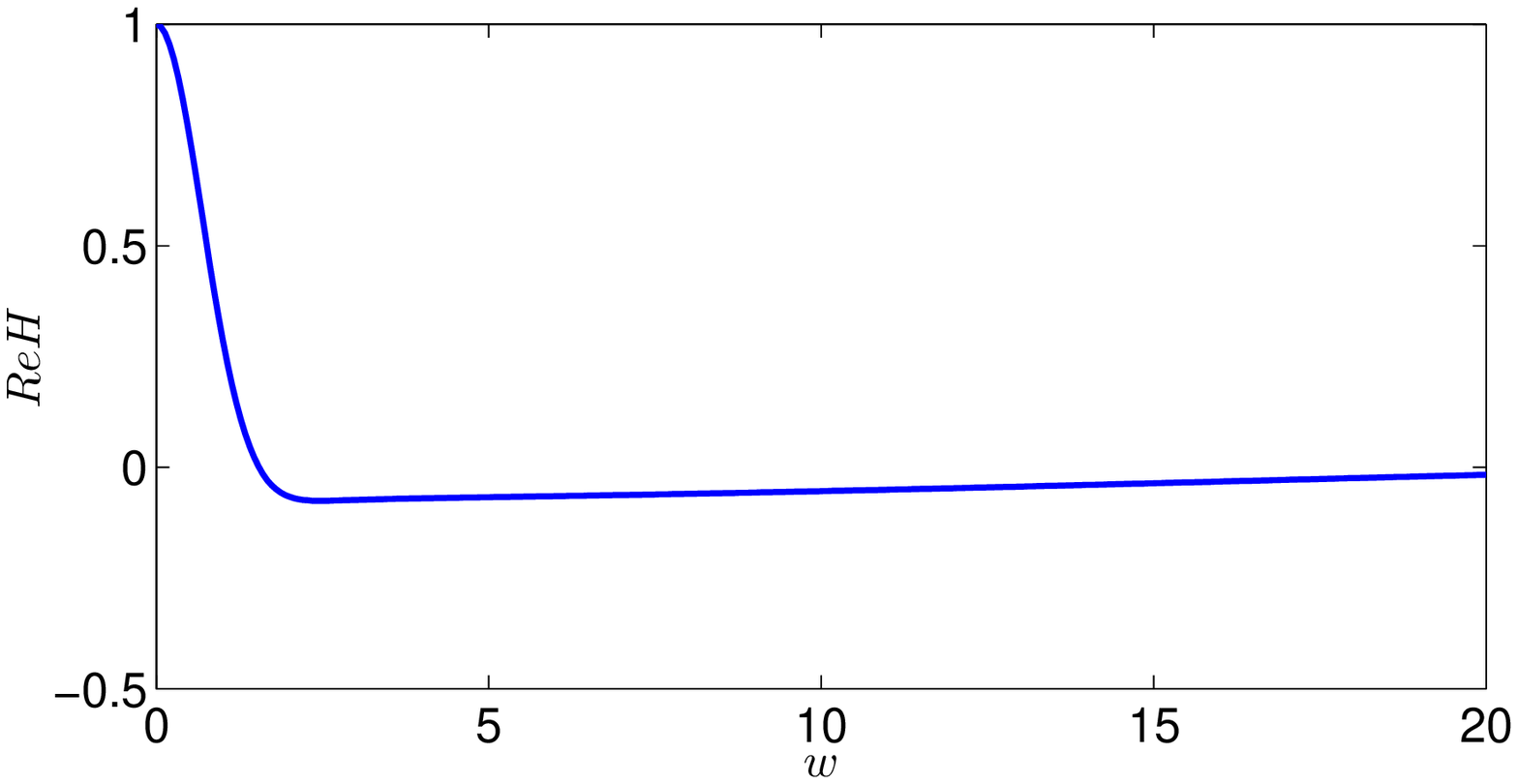}
\hspace{10pt}
\includegraphics[width=2.7in,angle=0]{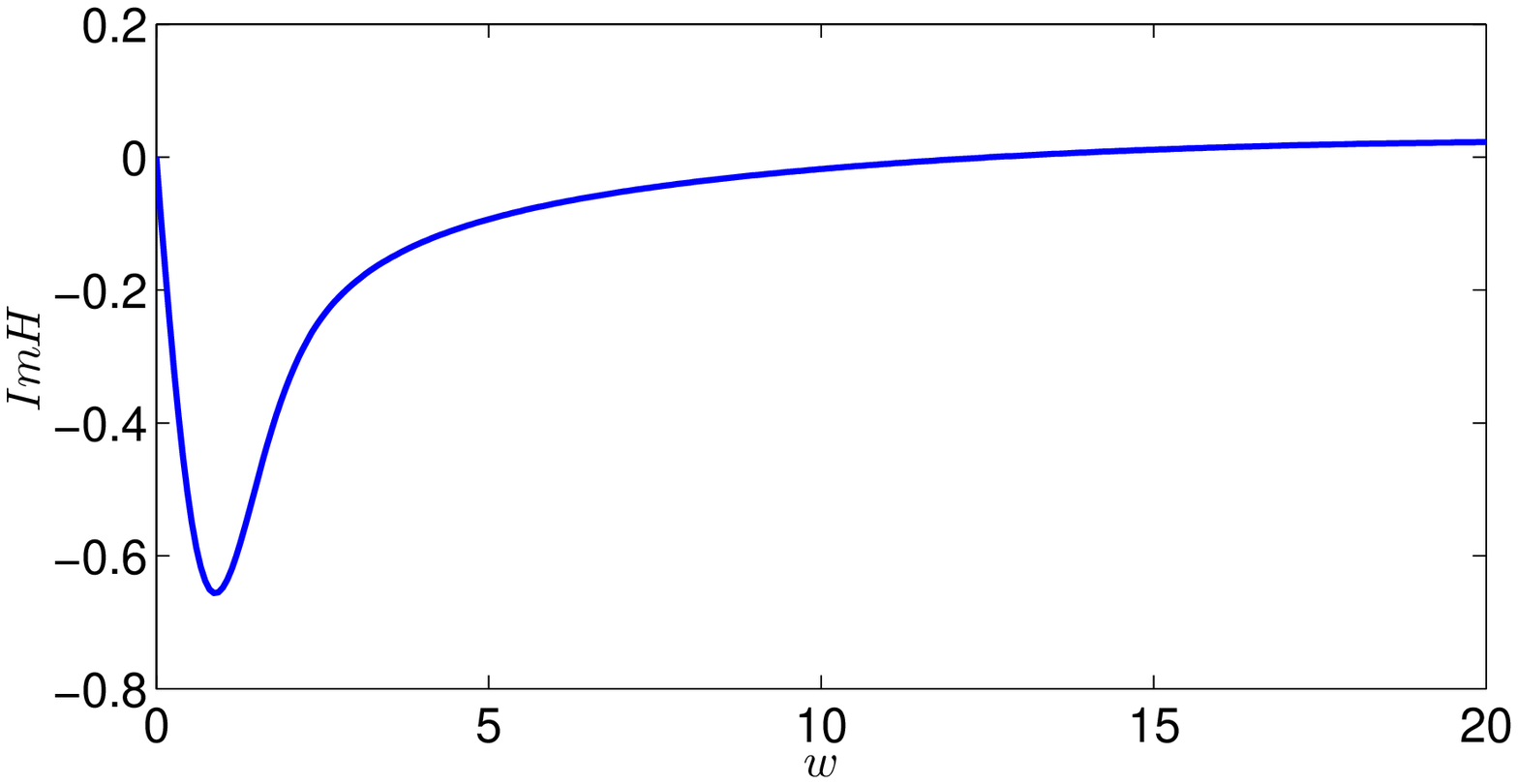}
\end{center}
\caption{$\Re H(w)$ and $\Im H(w)$ in the Dawson's integral example  using $N=300$ sample points.}
\label{FDawson1}
\end{figure}

\begin{figure} \begin{center}
\includegraphics[width=2.7in,angle=0]{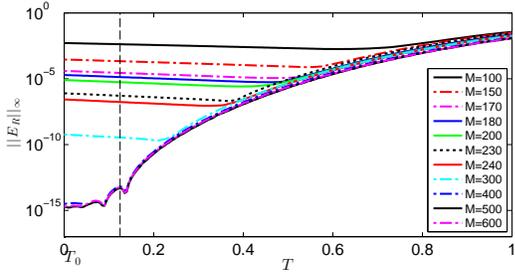}
\end{center}
\caption{$||E_R||_\infty$ in the Dawson's integral example as $M$ varies. Vertical dashed line indicates the time delay $T_0=0.125$.}
\label{FDawson2}
\end{figure}

We consider here another analytic example \cite{Knockaert_Dhaene_2008} modeled by the transfer function
\[
H(w)=\e^{-iwT_0}\tilde H(w)
%
\]
where 
\[
\tilde H(w)=\e^{-w^2}-\frac{2i}{\sqrt\pi} D(w),
\]
$D(w)$ is Dawson's integral
\[
D(w)=\e^{-w^2}\int_0^w \e^{t^2} dt = \frac{\sqrt\pi}{2}\e^{-w^2}\erfi(w)
\]
and $\erfi(w)$ is the imaginary error function.  Since \cite{Weideman_1995}\footnote{Please note that we use an opposite sign in the definition of the Hilbert transform than in \cite{Weideman_1995} and \cite{Knockaert_Dhaene_2008}.}
\[
{\mathcal H}[\Re \tilde H]={\mathcal H}(\e^{-w^2})= \frac{2}{\sqrt\pi} D(w)=-\Im \tilde H,
\]
function $\tilde H(w)$ is causal. Hence, the function $H(w)$ is  a causal function delayed with offset $T_0$. We use $T_0=0.125$ and sample $H(w)$ on $[0,w_{max}]$ with $w_{max}=20$ using various numbers of points ranging from $N=100$ to $600$. 
Real and imaginary parts of $H$ are shown in Fig. \ref{FDawson1}.
%
%
The evolution of $||E_R||_\infty$ for various $M$ is shown in Fig. \ref{FDawson2}, where it is clear that critical transition times approach $T_0$.
Constructing fitted quadratic error curves and extrapolating them to find their intersections with the horizontal line corresponding to  the error value $\xi$ produces a set of approximation of $T_0$, shown in Table \ref{T0approx_Dawson_table}.
%
%
\begin{table}[h]
\begin{center}
\begin{tabular}{|c|c||c|c|}
\hline
\rule{0cm}{10pt}
$M$ & $T_0$ estimate & $M$ & $T_0$ estimate \\[3pt]
\hline
\rule{0cm}{10pt}
  150 &     0.16735 & 240  & 0.12362 \\[3pt]
\hline
\rule{0cm}{10pt}  
  170 &  0.13116  &   300  &  0.12661 \\[3pt]
\hline
\rule{0cm}{10pt}  
180  &  0.12964 & 400  & 0.12233 \\[3pt] 
\hline
\rule{0cm}{10pt}    
200  &  0.12753 & 500  &  0.12578 \\
\hline
\rule{0cm}{10pt}    
230  &  0.12333 & 600  &  0.12414 \\
\hline
\end{tabular}
\end{center}
\caption{Approximations of $T_0$ in the Dawson's integral example using extrapolations of fitted error curves for various $M$. The exact value $T_0=0.125$, averaged value for $M$ ranging from 200 to 600 is $T_0^{(aver)}=0.12528$.}
\label{T0approx_Dawson_table}
\end{table} 
Averaging obtained approximations of $T_0$ for $M\geq 200$, once some convergence is established, gives $T_0^{(aver)}=0.12528$.

\begin{figure} \begin{center}
\includegraphics[width=2.7in,angle=0]{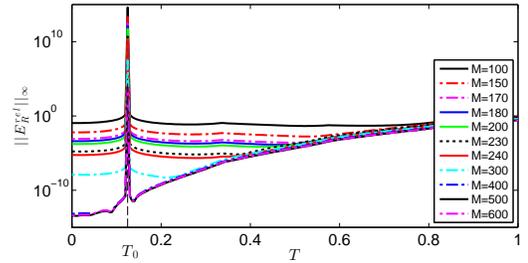}
\end{center}
\caption{Evolution of the relative error $||E_R^{rel}||_\infty$ in the Dawson's integral example as $M$ varies. Vertical dashed line indicates the time delay $T_0=0.125$.}
\label{FDawson3}
\end{figure}

It is interesting to note behavior of the relative error $E_R^{rel}$ in this example. The evolution of its $\infty$ norm is shown in Fig. \ref{FDawson3}. It is clear that all profiles even for small values of $M$ have a unique local maximum at $T=T_0$. $2$-norm has a similar behavior. Even though the behavior of the relative error $E_R^{rel}$  can be used to determine the time delay in this example, we did not find the same pronounced behavior in other examples we considered. At the same time, extrapolating fitted quadratic curves of $\infty$ norms of the absolute error $E_R$ was a robust approach in all examples we considered.
%

%

\subsection{Stripline Example} \label{Example_stripline_Wang}

\begin{figure} \begin{center}
\includegraphics[width=2.7in,angle=0]{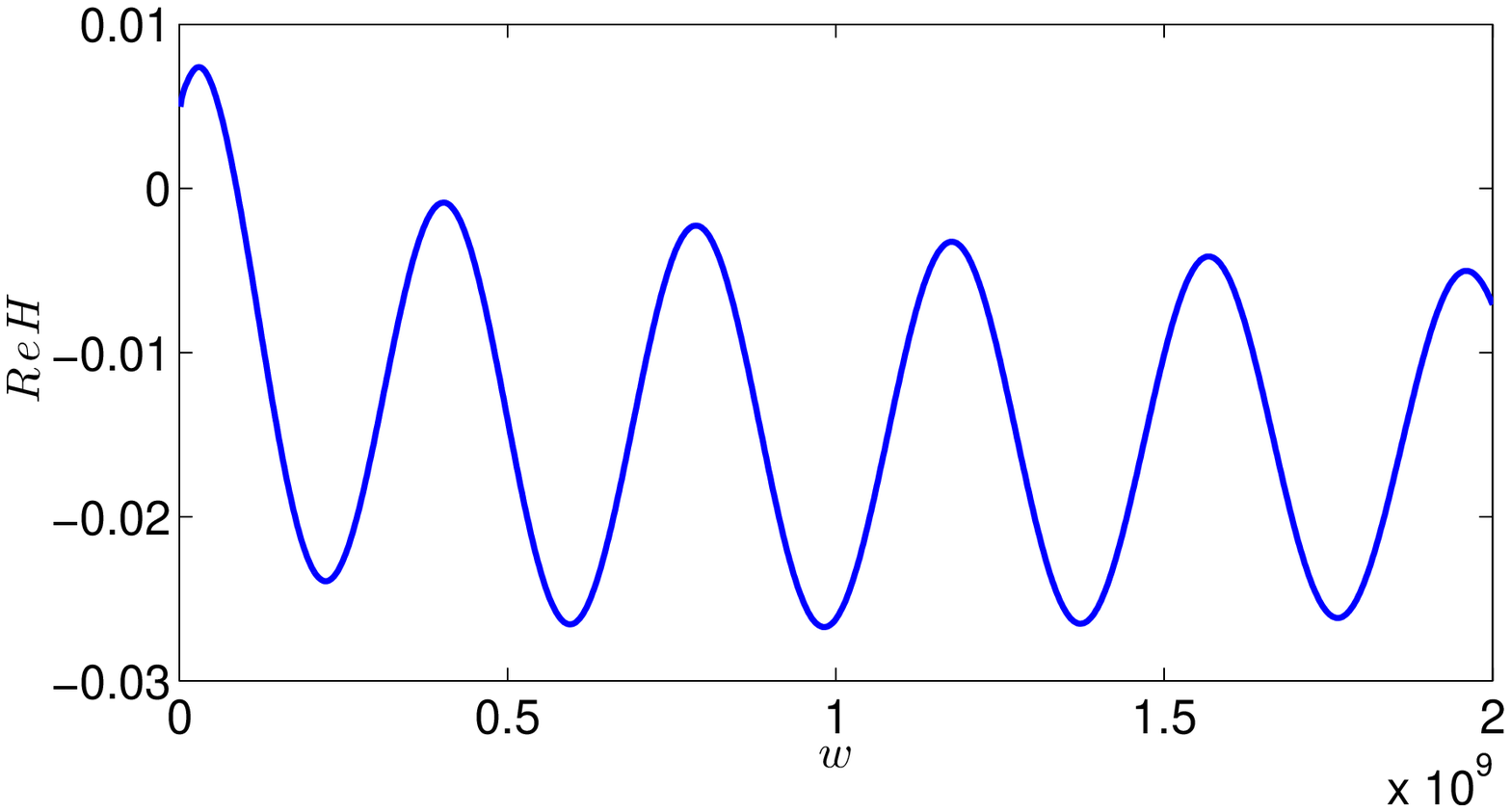}
\hspace{10pt}
\includegraphics[width=2.7in,angle=0]{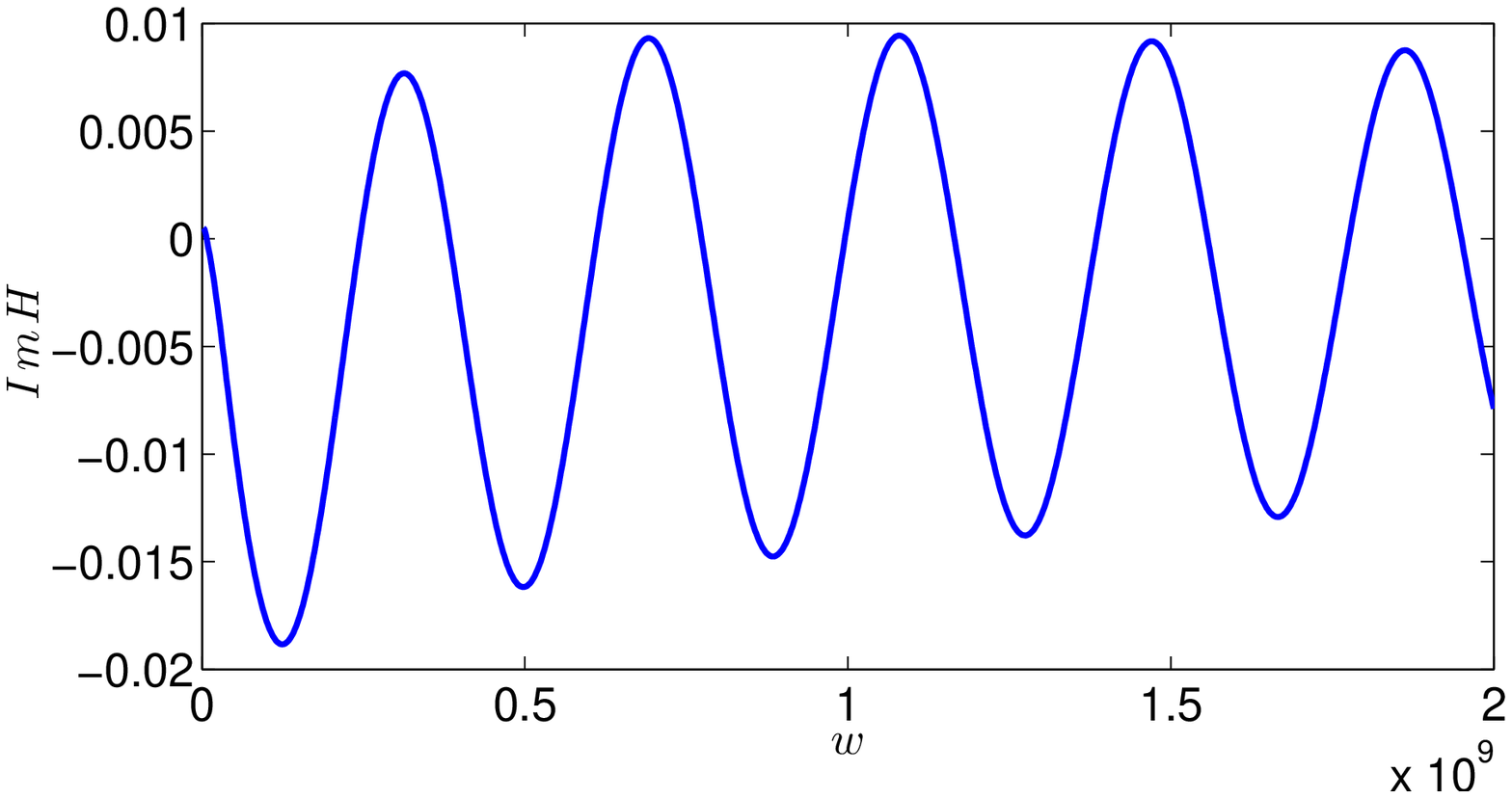}
\end{center}
\caption{$\Re H$ and $\Im H$ in the stripline  example with $N=1000$.}
\label{Wang0}
\end{figure}
\begin{figure} \begin{center}
\includegraphics[width=2.7in,angle=0]{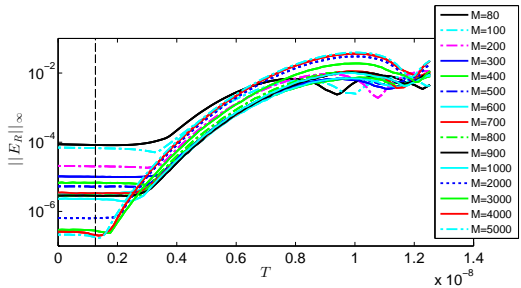}
\end{center}
\caption{$||E_R||_\infty$ in the stripline  example for various $M$. Vertical dashed line indicates the closed form microwave theory time delay approximation  $T_0=1.25809$\,ns.}
\label{Wang1}
\end{figure}

We simulated an asymmetric stripline  modeled in \cite{Wang_Drewniak_Fan_Knighten_Smith_Alexander_2002} 
with length $L=8$\,in, width $W=14$\,mils, distances from the trace to reference planes $H_1=10$\,mils, $H_2=20$\, mils, substrate dielectric  Megtron6-1035, Laminate with a dielectric constant $\epsilon_r= 3.45$ using a Cadence software tool with FEM full-wave field solver. The stripline was simulated on $[0,w_{max}]$ with $w_{max}=2$\,GHz. We analyzed element $H(w)=S_{11}(w)$ of the transfer matrix. The real and imaginary parts of $H$ are shown in Fig. \ref{Wang0}.
The evolution of $||E_R||_\infty$ for various $M$ is depicted in Fig. \ref{Wang1}. 
Even for high values of $M$, the error in causality does not go to machine precision and instead levels off around $10^{-6}$. This indicates that our finite element simulation results are accurate only within $10^{-7}-10^{-6}$. For causality characterization, this implies that data have noise/approximation errors with amplitude around $10^{-7}-10^{-6}$. Graphs of $||E_R||_\infty$ suggest that for $M\leq 2000$, the error is dominated by Fourier series approximation error, while for higher of $M$, the error is dominated by the noise/approximation errors from the finite element method.

\begin{figure} \begin{center}
\includegraphics[width=2.7in,angle=0]{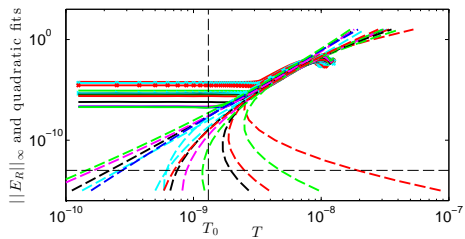}
\end{center}
\caption{Extrapolated quadratic curves based on initial  fitting range in the stripline example.}
\label{Wang2}
\end{figure}

\begin{figure} \begin{center}
\includegraphics[width=2.7in,angle=0]{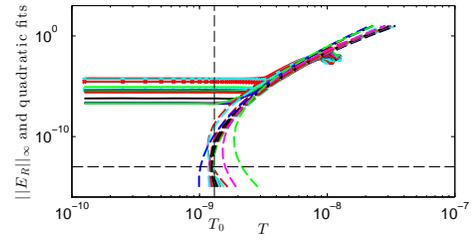}
\end{center}
\caption{Extrapolated quadratic curves based on more narrow  fitting range in the stripline example. The average time delay is $T_0^{(aver)}=1.2669$\,ns.}
\label{Wang3}
\end{figure}

In this example, the time delay was estimated using a closed form microwave theory approximation as $T_0=8 \times 0.0254/(c_0/ \sqrt{\epsilon _r})= 1.25809$\,ns, where $c_0=3\times 10^{8}$m/s is the speed of light, $0.0254$ is a conversion factor to convert from inches to meters. The error curves were fitted to quadratic curves as explained above. Because of relatively high errors in data, the fitted regions are not long enough. 
Besides, there is more nonlinear behavior of the error curves for high values of $T> T_0$. All this 
makes it difficult to estimate the time delay as shown in Fig. \ref{Wang2}. As can be seen, extrapolated quadratic curves do not focus at $T_0$ but instead spread out around $T_0$ similar to the four-pole example with an imposed noise considered in Section \ref{Example_four_pole}.
This problem can be corrected by decreasing the fitting range and going more away from transition regions. The results are shown in Fig. \ref{Wang3}. 
The approximations of $T_0$ are given in Table \ref{T0approx_stripline_Wang_table}. Averaging them for values of $M$ up to $3000$ produces
$T_0^{(aver)}=1.2669$\,ns, that agrees well with an analytically estimated time delay using a closed form formula. As in other examples, results with very high values of $M>3000$, that are more affected by noise and approximation errors in data, are less accurate.
\begin{table}[h]
\begin{center}
\begin{tabular}{|c|c||c|c|}
\hline
\rule{0cm}{10pt}
$M$ & $T_0$ estimate (in ns) & $M$ & $T_0$ estimate (in ns) \\[3pt]
\hline
\rule{0cm}{10pt}
  80 &     1.266 & 700  & 1.1987 \\[3pt]
\hline
\rule{0cm}{10pt}  
  100 &  1.2312  &   800  &  1.2179 \\[3pt]
\hline
\rule{0cm}{10pt}  
200  &  1.2553 & 900  & 1.2593 \\[3pt] 
\hline
\rule{0cm}{10pt}    
300  &  1.2797 & 1000  &  1.2205 \\
\hline
\rule{0cm}{10pt}    
400  &  1.2826 & 2000  &  1.246 \\
\hline
\rule{0cm}{10pt}    
500  &  1.2413 & 3000  &  1.5878 \\
\hline
\rule{0cm}{10pt}    
600  &  1.1833 & 4000  &  1.0168 \\
\hline
\end{tabular}
\end{center}
\caption{Approximations of $T_0$ in the stripline example using extrapolations of fitted quadratic error curves for various $M$. The closed form approximation of the time delay is $T_0=1.25809$\,ns, averaged value for $M$ ranging to up to $3000$ is $T_0^{(aver)}=1.2669$\,ns.}
\label{T0approx_stripline_Wang_table}
\end{table} 

\section{Conclusions} \label{conclusions}

We proposed a new method for time delay extraction from tabulated frequency responses. The approach uses
 the spectrally accurate causality enforcement technique constructed using SVD-based causal Fourier continuations, that was recently developed by the authors. The  time delay is incorporated to the causality characterization approach by introducing a linear varying phase factor to the system of equations that defines Fourier coefficients. Varying time until a threshold time, that depends on the maximum frequency at which the transfer function is available, results in the reconstruction error between the given data and their causal Fourier continuations to go from an almost constant small value to a rapidly growing function at some critical transition time. The critical transition times depend on resolution and approach the time delay as resolution increases. Several sets of frequency responses with increasing resolution can be used to establish convergence and get an approximation of the time delay. Alternatively, when only a limited number of samples is available, a growing portion of the error curve can be extrapolated to find an approximation of the time delay. The method is applicable to data that have noise or other approximation errors. A few sets of frequency responses can be used to improve the accuracy of time delay approximation by averaging the obtained results. The technique can be extended for multi-port and mixed  mode networks.
The performance of the method is demonstrated using several analytic and simulated examples, including data with noise, for which time delay is known exactly or can be estimated using other approaches.

\section*{Acknowledgment}

The authors are grateful to Dr. Linh V. Nguyen for valuable discussions on the Fourier transform. 
The work was supported by the Micron Foundation. 
The author L.L.B. would also  like to acknowledge the availability of computational resources made possible through the National Science Foundation Major Research Instrumentation Program, grant 1229766.

\section*{Appendix \\[7pt]
Error Analysis of Causality Characterization Method Based on Causal Fourier Continuations} \label{appendix}

In this section, we provide an upper bound of the reconstruction error between a given transfer function $H(x)$ and its causal Fourier continuation ${\mathcal C}(H)(x)$ in the presence of noise $\epsilon$ in data.

Denote by $\hat H_{M}$ any function of the form  
\begin{equation}\label{Fourier_series_arbitrary}
\hat H_{M}(x)=\sum_{k=1}^{M}  \hat\alpha_k \phi_k(x)
\end{equation}
where $\phi_k(x)=\e^{-\frac{2\pi i}{b}kx}$, $k=1,\ldots,M$. 

Let $A=U \Sigma V^*$ be the full SVD decomposition \cite{Trefethen_Bau_1997} of the matrix $A$ with entries $A_{kj}=\phi_k(x_j)$, $j=1,\ldots,N$, $k=1,\ldots, M$,
%
%
where $U$ is an $N\times N$  unitary matrix, 
$\Sigma$ is an $N\times M$ diagonal  matrix of singular values $\sigma_j$, $j=1,\ldots,p$, $p=\min({N,M})$,  $V$ is an $M \times M$ unitary matrix with entries $V_{kj}$,
 and $V^*$ denotes the complex conjugate transpose of $V$. 

The following result is true \cite{Barannyk_Aboutaleb_Elshabini_Barlow_IEEE}.

\begin{theorem} 

Consider a rescaled transfer function $H(x)$ defined by symmetry on $\Omega=[-0.5,-a]\cup[a, 0.5]$, where $a=0.5\frac{w_{min}}{w_{max}}$, whose values are available at points $x_j\in \Omega$, $j=1,\ldots,N$. Then the error in approximation of $H(x)$, that is known with some error $\varepsilon$, by 
its causal Fourier continuation ${\mathcal C}(H)(x)$ defined in (\ref{E3_0}) on a wider domain $\Omega^c=[-b/2,b/2]$, $b\geq 1$,
has the upper bound
\[
||H- {\mathcal C}(H+\varepsilon) ||_{L_2(\Omega)} \leq (1+\Lambda_2 \sqrt{N(M-K)}) 
\] 
\[
\times \left(|| H-\hat H_{M} || _{L_\infty(\Omega)} + ||\varepsilon||_{L_\infty(\Omega)} \right)
+\Lambda_1 \sqrt{K/b} ||\hat H_{M}||_{L_\infty(\Omega^c)}
\]
%
and holds for all functions of the form (\ref{Fourier_series_arbitrary}). Here
%
\[
\Lambda_1=\max_{j:\ \sigma_j<\xi}|| v_j(x) ||_{L_2(\Omega)}, \quad
\Lambda_2=\max_{j:\ \sigma_j>\xi}\frac{||v_j(x) ||_{L_2(\Omega)}}{\sigma_j},
\]
%
and  functions $v_j(x)=\sum_{k=1}^{M} V_{kj} \phi_k(x)$ are each an up to $M$ term causal Fourier series with coefficients given by the $j$th column of $V$; $K$ denotes the number of singular values that are discarded, i.e. the number of $j$ for which $\sigma_j<\xi$, where $\xi$ is the cut-off tolerance. 

\end{theorem}

It can be seen that constants $\Lambda_1$,  $\Lambda_2$  and $K$ depend only on the continuation parameters $N$, $M$, $b$ and $\xi$ as well as location of discrete points $x_j$, and not on the function $H$.  

For brevity, we can write the above error estimate 
 as
\[
||H- {\mathcal C}(H+\varepsilon) ||_{L_2(\Omega)} \leq \epsilon_F+\epsilon_n + \epsilon_T.
\]
Here
\[
\epsilon_F\equiv (1+\Lambda_2 \sqrt{N(M-K)}) || H-\hat H_{M} || _{L_\infty(\Omega)} 
\]
is the error due to a causal Fourier series approximation and it decays at least as fast as ${\mathcal O}(M^{-k+1})$, $k$ is the smoothness order of $H(x)$, which can be estimated numerically using reconstruction errors with different resolution (see \cite{Barannyk_Aboutaleb_Elshabini_Barlow_IEEE}).
%
\[
\epsilon_T=\Lambda_1 \sqrt{K/b} ||\hat H_{M}||_{L_\infty(\Omega^c)},
\]
%
that is the error due to the truncation of singular values. It is typically small and close to the cut-off value $\xi$.
\[
\epsilon_n=(1+\Lambda_2 \sqrt{N(M-K)}) ||\varepsilon||_{L_\infty(\Omega)} 
\]
is the error due to noise $\epsilon$ in data. Numerical experiments reveal that $\epsilon_n$ has the order of noise $\epsilon$ in data.




\bibliographystyle{acm}

\bibliography{references_time_delay}

\end{document}